\documentclass{article}

\usepackage[a4paper,top=2cm,bottom=2cm,left=3cm,right=3cm,marginparwidth=1.75cm]{geometry}
\usepackage{setspace}
\usepackage[english]{babel}
\usepackage{graphicx} % For including graphics
\usepackage{float} % To make figures appear exactly where you want
\usepackage{caption}
\usepackage{subcaption}
\usepackage{standalone} % Required for standalone TikZ files
\usepackage{tikz}
\usetikzlibrary{angles,quotes,calc,arrows.meta}
\usepackage{enumitem}
\usepackage{amsmath}
\usepackage{amssymb}
\usepackage{siunitx}
\usepackage[utf8]{inputenc}
\usepackage{csquotes}
\usepackage{booktabs}
\usepackage{longtable}
\usepackage{adjustbox}
\usepackage{array}
\usepackage{url}
\usepackage{titlesec}
\usepackage{indentfirst} % Indent the first paragraph of each section
\usepackage{authblk}
\usepackage{xcolor} % Load the xcolor package for color options

\usepackage[backend=biber,style=apa,natbib,bibencoding=utf8,sorting=nyt,hyperref=true,date=year,urldate=comp]{biblatex}
\usepackage{hyperref} % For clickable references
\usepackage{cleveref} % Enhanced references
\usepackage{xurl}  % Better automatic URL breaking
\hypersetup{
    colorlinks=true,        % Enable colored links
    linkcolor=blue,         % Color for internal links (figures, equations)
    filecolor=magenta,      % Color for file links
    urlcolor=cyan,          % Color for URLs
    citecolor=black,        % Make citations black but still clickable
    breaklinks=true
}
\newcommand{\figref}[1]{\textcolor{blue}{\hyperref[#1]{Fig.~\ref{#1}}}} % To make all text of Fig. 1 clickable

\titleformat{\subsection}
  {\mdseries\itshape\large} % Medium series, italic shape, and large font size
  {\thesubsection}{1em}{} % Numbering, spacing, and the section title itself

\title{A New Approach to Learn Trigonometry}
\author[1]{Marcia Ann Surya}
\author[2]{Yohanes Surya}
\affil[1]{Teachers College, Columbia University, New York, USA, mas2637@tc.columbia.edu}
\affil[2]{Gasing Academy, Jakarta, Indonesia, yohanes.surya@gasing.id}

\begin{document}
\maketitle

\begin{abstract}
We introduce the Primary Gasing Triangle, a right triangle with a hypotenuse of 1 unit, to define the primary trigonometric functions: sine and cosine. This triangle serves as the foundational element in a new approach to learning trigonometry, enabling us to derive the Derived Gasing Triangle, where the other four trigonometric functions (tangent, secant, cotangent, and cosecant) are defined. Using the Primary Gasing Triangle, we derive key trigonometric formulas, provide several proofs of the Pythagorean theorem using trigonometry, and solve various trigonometry problems. This approach makes learning trigonometry simpler, easier, and more intuitive.
\begin{center}
\textbf{Acknowledgment} \\[0.7em]
\end{center}
We would like to express our gratitude to Prof. Ben Koo for his insightful discussions, which greatly refined the concept of the Gasing Triangle in teaching trigonometry.
\end{abstract}

\section{Introduction}
\label{sec:Introduction}

Trigonometry is known as the study of angles and sides of a triangle. The name Trigonometry itself can be interpreted as triangle measurement. To facilitate the relationship between angles and triangles, six Trigonometric ratios or functions known as sine, cosine, tangent, cotangent, secant, and cosecant were introduced \autocite{gullberg_mathematics_1996, van_brummelen_2009}.

Initially, these six Trigonometric functions were defined as line segments \autocite{van_sickle_2011}. Later, this definition was changed to Trigonometric ratios, which represent the ratios of sides in a right triangle. This ratio definition is used in almost all Trigonometry textbooks \autocite{hall_knight_1928, gelfand_saul_2001}. However, this ratio-based definition does not provide a clear physical meaning and is difficult to apply for angles greater than 90 degrees. It also does not reflect the periodic nature of Trigonometric functions, preventing us from gaining a full understanding of various periodic phenomena such as waves, vibrations, light, sound, alternating current, or all other periodic phenomena, which is the primary importance of Trigonometry \autocite{simmons_1987}. 

The unclear physical meaning and applications of Trigonometry significantly affect students' motivation to learn it. Students memorize formulas without understanding their origins and quickly forget this subject after passing the test.

\begin{figure}[H]
    \centering
    % Include the TikZ figure from figure_1.tex
    \scalebox{1.1}{\input{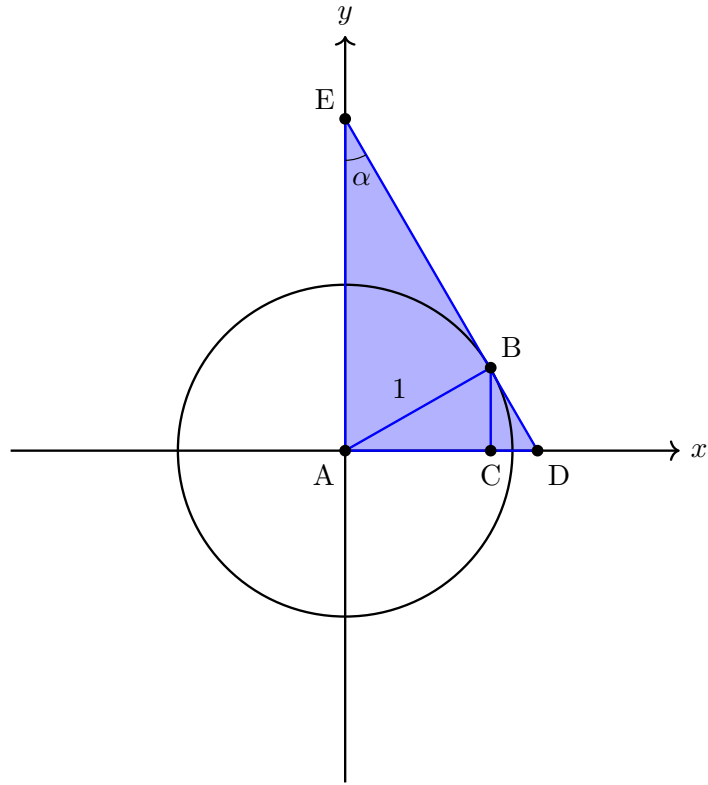}}
    \caption{Gasing Triangles}
    \label{Figure_1}
\end{figure}

In Fig.~\ref{Figure_1}, we introduce three similar triangles ABC, ADB, and EAB, which we call the \textbf{Gasing Triangle}, where one side of the triangle, AB, has a length of 1 unit. We name the triangle ABC the Primary Gasing Triangle with a hypotenuse length of 1 unit. The triangles ADB and EAB are introduced as Derived Gasing Triangles, the reasoning for which will be explained later.

We depict these Gasing Triangles in Cartesian coordinates and also illustrate a circle with a radius of 1 unit to facilitate the extension of the Trigonometric function definitions to obtuse angles.

The term Gasing in Segitiga Gasing is an acronym derived from the Indonesian phrase tiGA SIku pentiNG (meaning "three important right triangles"), referring to the three Gasing Triangles shown in Fig.~\ref{Figure_1}. Additionally, the word Gasing is often used to describe a learning method that is GAmpang (easy), aSIk (fun), and meNyenangkan (enjoyable) \autocite{gasing_academy}.

In this paper, we will demonstrate that the key element in learning Trigonometry lies in the Primary Gasing Triangle. Using this triangle, we will define and provide a physical meaning for the sine and cosine functions, which we refer to as the Primary Trigonometric Functions. Additionally, the Primary Gasing Triangle allows us to derive the other trigonometric functions, which we call Derived Trigonometric Functions, prove the Pythagorean theorem using Trigonometry, and derive essential Trigonometric formulas such as the sum and difference formulas, as well as the sine and cosine rules. We also extend the definition of the Primary Trigonometric Functions to obtuse angles. In \hyperref[appendix1]{Appendix 1}, we provide examples of solving Trigonometry problems using the Primary Gasing Triangle, while in \hyperref[appendix2]{Appendix 2}, we offer several proofs of the Pythagorean theorem using the Primary Trigonometric Functions, which were previously thought to be impossible \autocite{loomis_1968}.

\break
\section{Defining Trigonometric Functions}
\label{sec:Defining Trigonometric Functions}

\begin{figure}[H]
    \centering
    \scalebox{1.6}{\input{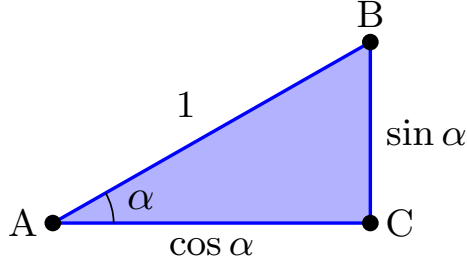}}
    \caption{Primary Gasing Triangle}
    \label{Figure_2}
\end{figure}

First, we draw the Primary Gasing Triangle ABC from Fig.~\ref{Figure_1}, as shown in Fig.~\ref{Figure_2}. Here, the length of side AB is 1 unit, the length of side BC is defined as sin $\alpha$, and the length of side AC is defined as cos $\alpha$, with $\alpha$ being an acute angle ($0\leq \alpha \leq 90^{\circ}$). It is clear that we define sin $\alpha$ as the length of the vertical side (the side opposite angle $\alpha$) and cos $\alpha$ as the length of the horizontal side (the side adjacent to angle $\alpha$) in the Primary Gasing Triangle ABC.
%~ makes them as one word so it won't break the two in two lines

For obtuse angles, the definitions of sine and cosine are extended as the projections of a unit vector onto the y-axis and x-axis in a Cartesian coordinate system. This extension remains consistent with the previous definition. The unit vector approach we use here differs from the unit circle approach \autocite{sullivan_2016, stewart_redlin_watson_2016}. While the unit circle approach can explain Trigonometric functions for obtuse angles, it lacks a clear physical meaning, particularly in applications related to physics or engineering \autocite{simmons_1987}.

The other four trigonometric functions—secant, cosecant, tangent, and cotangent—are defined as the sides of two triangles, ADB and EAB, as shown in Fig.~\ref{Figure_1}. To calculate the unknown sides of triangles ADB and EAB, we use the Primary Gasing Triangle ABC as a reference.

\begin{figure}[H]
    \centering
    \scalebox{1.4}{\input{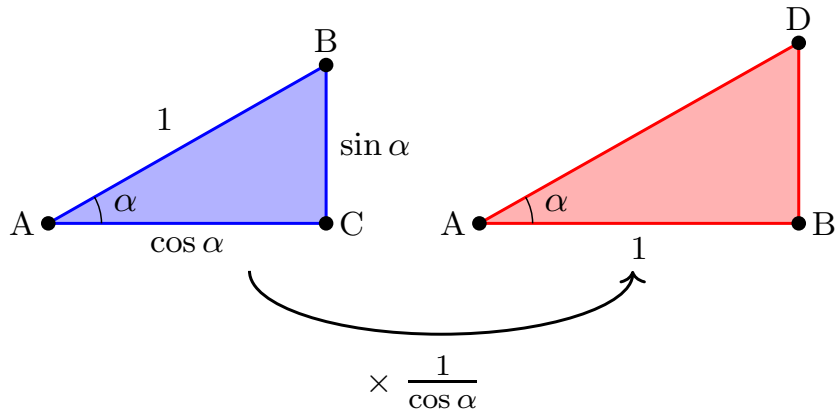}}
    \renewcommand{\thefigure}{3a} % Temporarily redefine figure numbering
    \captionsetup{labelformat=empty} % Removes default label format
    \caption{Figure 3a: Transformation of Primary Gasing Triangle ABC into Triangle ADB with a scale factor of $\frac{1}{\cos \alpha}$}
    \label{Figure_3a}
    \renewcommand{\thefigure}{\arabic{figure}} % Restore standard numbering
\end{figure}

In Fig.~\ref{Figure_3a}, the Primary Gasing Triangle ABC and the similar triangle ADB are depicted.

For similar triangles, the corresponding sides are connected by the same scaling factor. Since side AC = $\cos \alpha$ is connected to side AB = 1 unit by a scaling factor of $\frac{1}{\cos \alpha}$ , 
sides AD and BD can be calculated by multiplying this scaling factor by sides AB and CB, respectively. 

We obtain: 
\\
AD = AB $\times$ scale factor = $1 \times \frac{1}{\cos \alpha}$, hence AD = $\frac{1}{\cos \alpha}$
\\
BD = CB $\times$ scale factor = $\sin \alpha \times \frac{1}{\cos \alpha}$, hence BD = $\frac{\sin \alpha}{\cos \alpha}$

\begin{figure}[H]
    \centering
    \scalebox{1.6}{\input{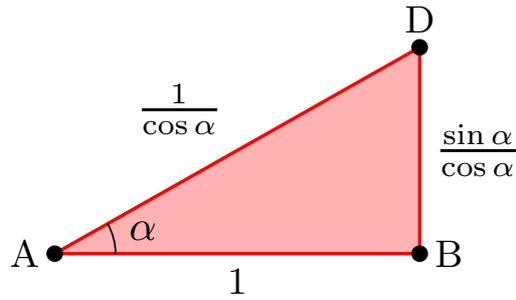}}
    \renewcommand{\thefigure}{3b} % Temporarily redefine figure numbering
    \captionsetup{labelformat=empty} % Removes default label format
    \caption{Figure 3b: Triangle ADB with side lengths derived from the transformation of Triangle ABC with a scale factor of $\frac{1}{\cos \alpha}$}
    \label{Figure_3b}
    \renewcommand{\thefigure}{\arabic{figure}} % Restore standard numbering
\end{figure}

This result is illustrated in Fig.~\ref{Figure_3b}.

Similarly, we can use the same method to calculate the unknown sides of triangle EAB.

\begin{figure}[H]
    \centering
    \scalebox{1.4}{\input{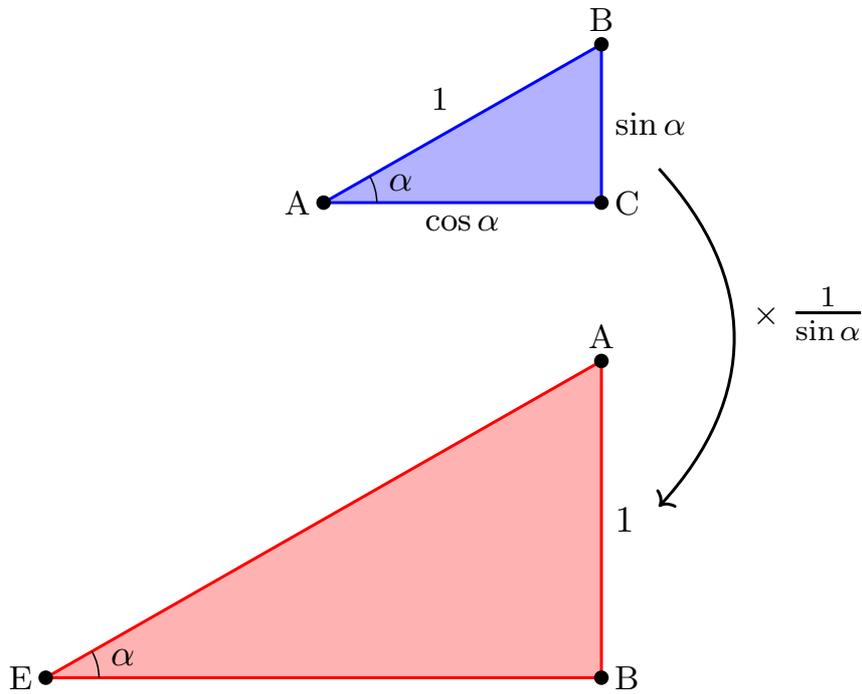}}
    \renewcommand{\thefigure}{4a} % Temporarily redefine figure numbering
    \captionsetup{labelformat=empty} % Removes default label format
    \caption{Figure 4a: Transformation of Primary Gasing Triangle ABC into Triangle EAB with a scale factor of $\frac{1}{\sin \alpha}$}
    \label{Figure_4a}
    \renewcommand{\thefigure}{\arabic{figure}} % Restore standard numbering
\end{figure}

In Fig.~\ref{Figure_4a}, we observe that the scaling factor is $\frac{1}{\sin \alpha}$.
Since triangles ABC and EAB are similar, the lengths of sides EA and EB can be obtained by multiplying the lengths of sides AB and AC by this scaling factor.
We obtain:
\\
EA = AB $\times$ scale factor = $1 \times \frac{1}{\sin \alpha}$, hence EA = $\frac{1}{\sin \alpha}$. 
\\
EB = AC $\times$ scale factor = $\cos \alpha \times \frac{1}{\sin \alpha}$, hence EB = $\frac{\cos \alpha}{\sin \alpha}$

\begin{figure}[H]
    \centering
    \scalebox{1.4}{\input{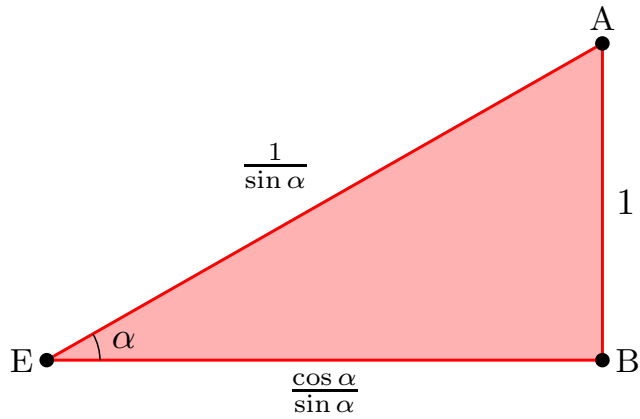}}
    \renewcommand{\thefigure}{4b} % Temporarily redefine figure numbering
    \captionsetup{labelformat=empty} % Removes default label format
    \caption{Figure 4b: Triangle EAB with side lengths derived from the transformation of Triangle ABC with a scale factor of $\frac{1}{\sin \alpha}$}
    \label{Figure_4b}
    \renewcommand{\thefigure}{\arabic{figure}} % Restore standard numbering
\end{figure}

This result is shown in Fig.~\ref{Figure_4b}.

Now, we define the side lengths as follows: AD = sec $\alpha$, BD = tan $\alpha$, EB = cot $\alpha$, and EA = csc $\alpha$. This allows us to establish the relationships between these four trigonometric functions and the Primary Trigonometric Functions, sine and cosine, that we defined earlier.
These relationships are expressed in equations (1) through (4):

\begin{align}
\tan \alpha &= \frac{\sin \alpha}{\cos \alpha} \tag{1} \label{eq:1} \\[0.8em]
\sec \alpha &= \frac{1}{\cos \alpha} \tag{2} \label{eq:2} \\[0.8em]
\cot \alpha &= \frac{\cos \alpha}{\sin \alpha} \tag{3} \label{eq:3} \\[0.8em]
\csc \alpha &= \frac{1}{\sin \alpha} \tag{4} \label{eq:4}
\end{align}

Fig.~\ref{Figure_5}a--c and Fig.~\ref{Figure_6} illustrate the six trigonometric functions, represented as the sides of the three Gasing Triangles.

% Reset figure counter to 4 so that Figure 5 is correctly numbered
\setcounter{figure}{4}
\begin{figure}[H]
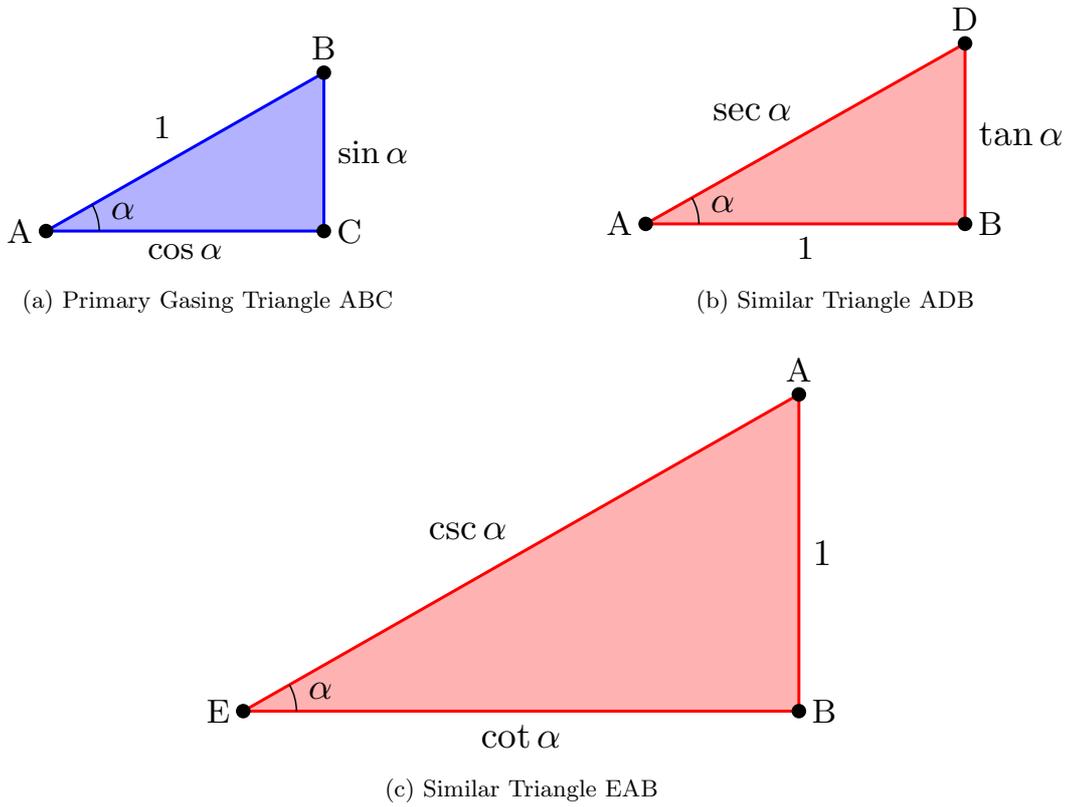

    \centering

    % Top row: Figure 5a and 5b
    \begin{subfigure}{0.45\textwidth}
        \centering
        \scalebox{1.4}{\input{figures/Figure_5a}}
        \caption{Primary Gasing Triangle ABC}
        \label{Figure_5a}
    \end{subfigure}
    \hfill
    \begin{subfigure}{0.45\textwidth}
        \centering
        \scalebox{1.4}{\input{figures/Figure_5b}}
        \caption{Similar Triangle ADB}
        \label{Figure_5b}
    \end{subfigure}

    % Bottom row: Figure 5c
    \vspace{1em} % Adjust spacing between rows
    \begin{subfigure}{0.9\textwidth}
        \centering
        \scalebox{1.4}{\input{figures/Figure_5c}}
        \caption{Similar Triangle EAB}
        \label{Figure_5c}
    \end{subfigure}

    \caption{The six trigonometric functions}
    \label{Figure_5}
\end{figure}

\begin{figure}[H]
    \centering
    \scalebox{1}{\input{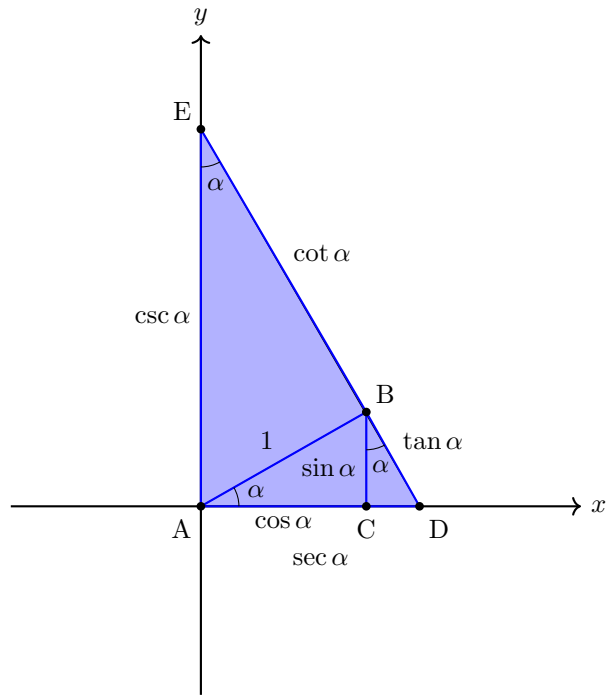}}
    \caption{The six trigonometric functions combined into a single diagram}
    \label{Figure_6}
\end{figure}

From the explanation above, the side lengths of the Gasing Triangles ADB and EAB in Figs.~\ref{Figure_5b} and ~\ref{Figure_5c} are derived from the side lengths of the Primary Gasing Triangle ABC in Fig.~\ref{Figure_5a}. Therefore, we refer to triangle ABC as the Primary Gasing Triangle, while ADB and EAB are referred to as Derived Gasing Triangles.

We also observe that the trigonometric functions secant, tangent, cosecant, and cotangent can be expressed in terms of the sine and cosine functions. Therefore, sine and cosine are called Primary Trigonometric Functions, while the other four functions are referred to as Derived Trigonometric Functions.

For all subsequent trigonometric calculations or derivations of trigonometric formulas in the remainder of this paper, we will rely solely on the Primary Gasing Triangle and the appropriate scaling factors.

\break
\section{Proving the Pythagorean Theorem with Trigonometry}
\label{sec:Proving the Pythagorean Theorem with Trigonometry}
Elisha Scott Loomis claimed that it is impossible to prove the Pythagorean Theorem using trigonometry \autocite{loomis_1968}. However, more recently, there have been several proofs of the Pythagorean Theorem using trigonometry \autocite{sloman_2023, shirali_2023, de_villiers_2023}. 
In this section, we will provide a trigonometric proof of the Pythagorean Theorem using the Primary Gasing Triangle. 
To do so, we will first calculate the side lengths of all the triangles depicted in Fig.~\ref{Figure_1}.

The lengths of sides $AD$ and $DE$ can be determined by the similarity between triangle $\triangle EDA$ and the Primary Gasing Triangle $\triangle ABC$, using a scaling factor of $\frac{1}{\sin \alpha \cos \alpha}$.  
We find that the length of $AD = \frac{1}{\cos \alpha}$ and $DE = \frac{1}{\sin \alpha \cos \alpha}$.  

For the sides $BD$ and $CD$, we use the similarity between the Primary Gasing Triangle $\triangle ABC$ and triangle $\triangle BDC$, applying a scaling factor of $\frac{\sin \alpha}{\cos \alpha}$.  
This gives us the lengths $BD = \frac{\sin \alpha}{\cos \alpha}$ and $CD = \frac{\sin^2 \alpha}{\cos \alpha}$.

\begin{figure}[H]
    \centering
    \scalebox{1}{\input{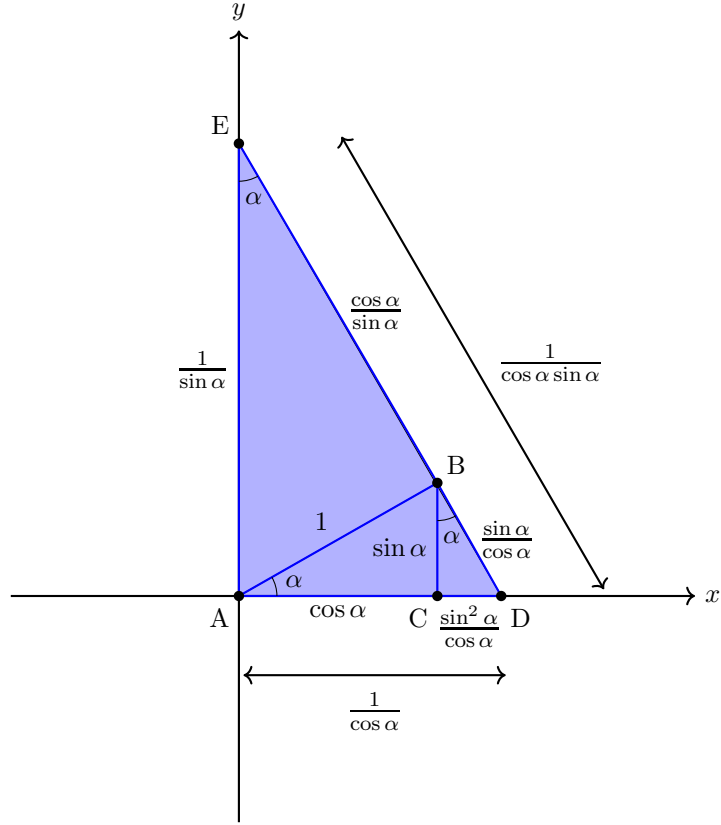}}
    \caption{Derived diagram from Figure 1 using ratio}
    \label{Figure_7}
\end{figure}

The complete results for all the side lengths of the triangles in Fig.~\ref{Figure_1} are shown in Fig.~\ref{Figure_7}. It is important to note that these calculations are based solely on the use of side ratios (scaling factors) and the similarity of right triangles.

From Fig.~\ref{Figure_7}, we observe the following:
\begin{align}
AD &= AC + CD \tag{5} \label{eq:5} \\
\frac{1}{\cos \alpha} &= \cos \alpha + \frac{\sin^2 \alpha}{\cos \alpha} \tag{6} \label{eq:6}
\end{align} % &= is specific for align only

By equalizing the denominators, we obtain
\begin{equation}
\frac{1}{\cos \alpha} = \frac{\cos^2 \alpha}{\cos \alpha} + \frac{\sin^2 \alpha}{\cos \alpha}, \tag{7} \label{eq:7}
\end{equation}
which simplifies to
\begin{equation}
1 = \cos^2 \alpha + \sin^2 \alpha. \tag{8} \label{eq:8}
\end{equation}

This result is a fundamental trigonometric identity. Notably, we obtained this result using only the similarity and the ratios of the sides of the triangles.

Using the Primary Gasing Triangle \( \triangle ABC \) (Fig.~\ref{Figure_5a}), we see that Eq.~\eqref{eq:8} is, in fact, the Pythagorean Theorem:
\begin{equation}
AB^2 = AC^2 + CB^2 \tag{9} \label{eq:9}
\end{equation}

Thus, we see that by proving the trigonometric identity, we simultaneously prove the Pythagorean Theorem. Therefore, there is no circular reasoning involved.

We can also prove this trigonometric identity (the Pythagorean Theorem) using Fig.~\ref{Figure_7}, where $DE = DB + BE$.

Explicitly, this can be written as
\begin{equation}
\frac{1}{\sin \alpha \cos \alpha} = \frac{\sin \alpha}{\cos \alpha} + \frac{\cos \alpha}{\sin \alpha}. \tag{10} \label{eq:10}
\end{equation}

By equalizing the denominators, we obtain
\begin{equation}
\frac{1}{\sin \alpha \cos \alpha} = \frac{\sin^2 \alpha}{\cos \alpha \sin \alpha} + \frac{\cos^2 \alpha}{\sin \alpha \cos \alpha}. \tag{11} \label{eq:11}
\end{equation}

Eq.~\eqref{eq:11} simplifies to $1 = \cos^2 \alpha + \sin^2 \alpha$, just as we expected.

Here, we see that the Pythagorean Theorem can indeed be proven using trigonometry without any concern for circular reasoning.

In \hyperref[appendix2]{Appendix 2}, we provide additional proofs of the Pythagorean Theorem using trigonometry. We observe that all valid proofs of the Pythagorean Theorem, whether geometric or algebraic \autocite{loomis_1968}, can be demonstrated through trigonometric proofs using the Primary Gasing Triangle and the appropriate scaling factors.

With the Pythagorean Theorem now proven, we can easily demonstrate the following identities using Figs.~\ref{Figure_5b} and \ref{Figure_5c}:
\begin{align}
\sec^2 \alpha &= 1 + \tan^2 \alpha \tag{12} \label{eq:12} \\
\csc^2 \alpha &= 1 + \cot^2 \alpha \tag{13} \label{eq:13}
\end{align}

\break
\section{Sum and Difference Formulas}
\label{sec:Sum and Difference Formulas}
In this section, we will derive the formulas for $\cos\,(\alpha \pm \beta)$ and $\sin\,(\alpha \pm \beta)$ using the Primary Gasing Triangle.

In Fig.~\ref{Figure_8a}, we observe the Primary Gasing Triangle $\triangle ABC$ with $\angle BAC = \alpha$. The hypotenuse $AB$ is 1 unit long. The side opposite angle $\alpha$ is $BC = \sin \alpha$, and the side adjacent to angle $\alpha$ is $AC = \cos \alpha$. We then construct a right triangle $\triangle ACF$, with $F$ as the right angle and $\angle CAF = \beta$.

\begin{figure}[H]
    \centering

    \begin{subfigure}{0.45\textwidth}
        \centering
        \scalebox{1.4}{\input{figures/Figure_8a}}
        \caption{Constructing Triangle ACF with side AC as one of the sides of the Primary Gasing Triangle ABC}
        \label{Figure_8a}
    \end{subfigure}
    \hfill
    \begin{subfigure}{0.45\textwidth}
        \centering
        \scalebox{1.4}{\input{figures/Figure_8b}}
        \caption{Gasing Triangle \( A_1C_1F_1 \) with angle \( \beta \)}
        \label{Figure_8b}
    \end{subfigure}

    \caption{}
    \label{Figure_8}
\end{figure}

To calculate the lengths of $AF$ and $CF$, we use the Primary Gasing Triangle $\triangle A_1C_1F_1$ (Fig.~\ref{Figure_8b}), which is similar to $\triangle ACF$, with the hypotenuse $A_1C_1 = 1$ unit. The side opposite angle $\beta$ is $F_1C_1 = \sin \beta$, and the side adjacent to angle $\beta$ is $A_1F_1 = \cos \beta$.

The Primary Gasing Triangle $\triangle A_1C_1F_1$ and $\triangle ACF$ are related by a scaling factor of $\cos \alpha$. Therefore, the lengths of sides $FC$ and $AF$ can be obtained by multiplying this scaling factor by the lengths $F_1C_1$ and $A_1F_1$. The results are
\[
FC = \cos \alpha \sin \beta \quad \text{and} \quad AF = \cos \alpha \cos \beta.
\]

\begin{figure}[H]
    \centering
    \scalebox{1.4}{\input{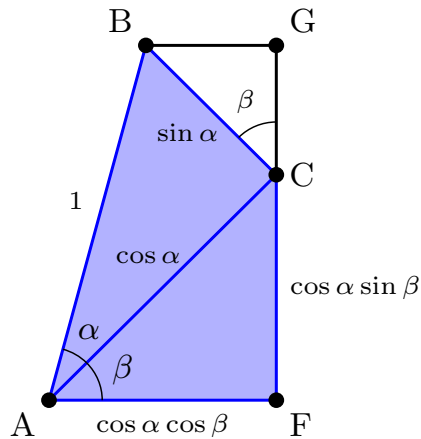}}
    \renewcommand{\thefigure}{8c} % Temporarily redefine figure numbering
    \captionsetup{labelformat=empty} % Removes default label format
    \caption{Figure 8c: Side lengths of Triangle ACF}
    \label{Figure_8c}
    \renewcommand{\thefigure}{\arabic{figure}} % Restore standard numbering
\end{figure}

Next, we construct a right triangle $\triangle CBG$ with a perpendicular at $G$, where $\angle BCG = \beta$, as shown in Fig.~\ref{Figure_8c}.

Since $\triangle CBG$ is similar to the Primary Gasing Triangle $\triangle A_1C_1F_1$ with a scaling factor of $\sin \alpha$, we can calculate the lengths of $CG$ and $GB$ by multiplying this scaling factor by the lengths $A_1F_1$ and $F_1C_1$. The results are
\[
CG = \sin \alpha \cos \beta \quad \text{and} \quad GB = \sin \alpha \sin \beta.
\]

\begin{figure}[H]
    \centering
    \scalebox{1.4}{\input{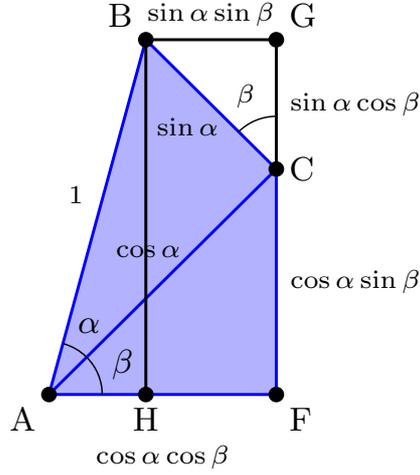}}
    \renewcommand{\thefigure}{8d} % Temporarily redefine figure numbering
    \captionsetup{labelformat=empty} % Removes default label format
    \caption{Figure 8d: Triangles \( ABC \), \( ACF \), and \( CBG \) with side lengths for calculating the Sum Formula}
    \label{Figure_8d}
    \renewcommand{\thefigure}{\arabic{figure}} % Restore standard numbering
\end{figure}

Next, we draw line $HB$. Notice that $\triangle ABH$ is a Primary Gasing Triangle with its hypotenuse $AB = 1$ unit and $\angle BAH = \alpha + \beta$, as shown in Fig.~\ref{Figure_8d}.

Here, $AH = \cos\,(\alpha + \beta)$ and $HB = \sin\,(\alpha + \beta)$. From Fig.~\ref{Figure_8d}, we can see that
\[
HB = FC + CG,
\]
or
\begin{equation}
\sin(\alpha + \beta) = \cos \alpha \sin \beta + \sin \alpha \cos \beta. \tag{14} \label{eq:14}
\end{equation}

By setting $\alpha = \beta$, we obtain
\begin{equation}
\sin(2\alpha) = 2 \sin \alpha \cos \alpha. \tag{15} \label{eq:15}
\end{equation}

Also, from Fig.~\ref{Figure_8d}, we have
\[
AH = AF - HF,
\]
or
\begin{equation}
\cos(\alpha + \beta) = \cos \alpha \cos \beta - \sin \alpha \sin \beta. \tag{16} \label{eq:16}
\end{equation}

By setting $\alpha = \beta$, we obtain
\begin{equation}
\cos(2\alpha) = \cos^2 \alpha - \sin^2 \alpha. \tag{17} \label{eq:17}
\end{equation}

Next, by applying the identity $\cos^2 \alpha + \sin^2 \alpha = 1$ to Eq.~\eqref{eq:17}, we obtain the following:
\begin{align}
\cos(2\alpha) &= 2 \cos^2 \alpha - 1 \tag{18} \label{eq:18} \\
\cos(2\alpha) &= 1 - 2 \sin^2 \alpha \tag{19} \label{eq:19}
\end{align}

Equations~\eqref{eq:15}, \eqref{eq:17}, \eqref{eq:18}, and \eqref{eq:19} are known as the double-angle identities.

To derive the formulas for $\sin\,(\alpha - \beta)$ and $\cos\,(\alpha - \beta)$, we can use Eqs.~\eqref{eq:14} and \eqref{eq:16} by substituting $\beta$ with $-\beta$.

\begin{figure}[H]
    \centering
    \scalebox{1.6}{\input{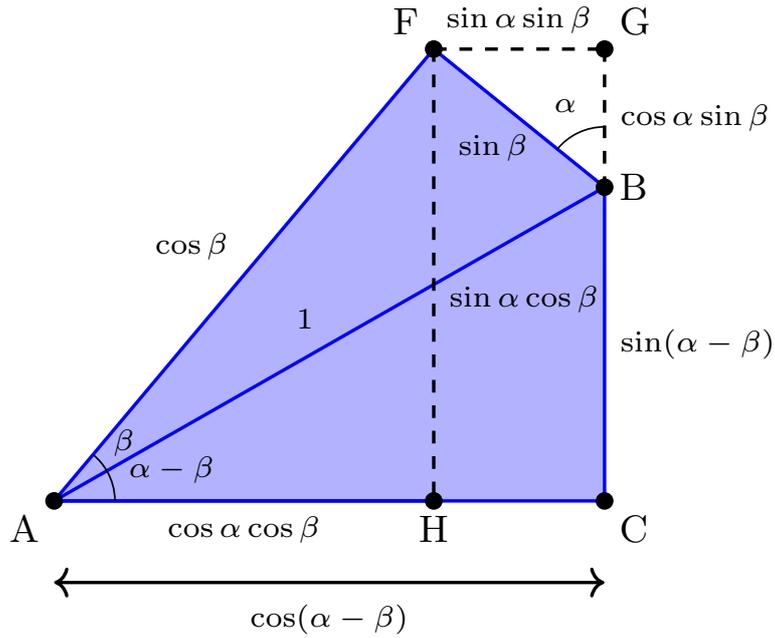}}
    \renewcommand{\thefigure}{8e} % Temporarily redefine figure numbering
    \captionsetup{labelformat=empty} % Removes default label format
    \caption{Figure 8e: Triangles \( ABC \), \( ABF \), and \( BFG \) with side lengths for calculating the Difference Formula}
    \label{Figure_8e}
    \renewcommand{\thefigure}{\arabic{figure}} % Restore standard numbering
\end{figure}

Alternatively, we can use a similar approach as above, but with the aid of Fig.~\ref{Figure_8e} (we leave the proof of these side lengths as an exercise for the readers).

In the figure, it is clear that
\[
CB = HF - BG,
\]
so
\begin{equation}
\sin\,(\alpha - \beta) = \sin \alpha \cos \beta - \cos \alpha \sin \beta. \tag{20} \label{eq:20}
\end{equation}

Similarly,
\[
AC = AH + HC,
\]
so
\begin{equation}
\cos\,(\alpha - \beta) = \cos \alpha \cos \beta + \sin \alpha \sin \beta. \tag{21} \label{eq:21}
\end{equation}

We can see how beautifully the sum and difference identities for cosine and sine are derived using the Primary Gasing Triangle and the appropriate scaling factors.

\break
\section{Sine and Cosine Rules}
\label{sec:Sine and Cosine Rules}
The sine and cosine rules can be derived using the Primary Gasing Triangle with the following steps:

First, we construct two Primary Gasing Triangles, $\triangle ABD$ and $\triangle CFE$, as shown in Figs.~\ref{Figure_9a} and~\ref{Figure_9b}.

\setcounter{figure}{8} % Reset figure counter to 8 so that Figure 9 is correctly numbered
\begin{figure}[H]
    \centering
    \begin{subfigure}{0.45\textwidth}
        \centering
        \scalebox{1.4}{\input{figures/Figure_9a}}
        \caption{Gasing Triangle with angle $\alpha$}
        \label{Figure_9a}
    \end{subfigure}
    \hfill
    \begin{subfigure}{0.45\textwidth}
        \centering
        \scalebox{1.4}{\input{figures/Figure_9b}}
        \caption{Gasing Triangle with angle $\gamma$}
        \label{Figure_9b}
    \end{subfigure}
    \caption{}
    \label{Figure_9}
\end{figure}

\begin{figure}[H]
    \centering
    \scalebox{1.4}{\input{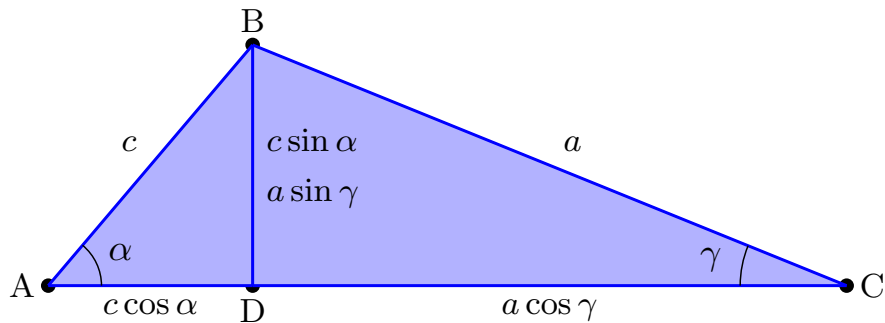}}
    \renewcommand{\thefigure}{9c} % Temporarily redefine figure numbering
    \captionsetup{labelformat=empty} % Removes default label format
    \caption{Figure 9c: Triangles \( ABD \) and \( CBD \) for deriving the Sine and Cosine Formulas}
    \label{Figure_9c}
    \renewcommand{\thefigure}{\arabic{figure}} % Restore standard numbering
\end{figure}

Next, we multiply the sides of $\triangle ABD$ by the scaling factor $a$ and the sides of $\triangle CFE$ by the scaling factor $c$, so that the new length $DB$ is equal to $EF$. Now we have a new triangle $\triangle ABC$, which is not a right triangle (see Fig.~\ref{Figure_9c}).

From Fig.~\ref{Figure_9c}, it can be seen that the side $DB$ satisfies
\begin{equation}
c \sin \alpha = a \sin \gamma. \tag{22} \label{eq:22}
\end{equation}

From this, we obtain the sine rule:
\begin{equation}
\frac{a}{\sin \alpha} = \frac{c}{\sin \gamma} \tag{23} \label{eq:23}
\end{equation}

Now, let's define the length $AC = b$. Since $AC = AD + DC$, we have
\[
b = c \cos \alpha + a \cos \gamma.
\]

In $\triangle CBD$, we use the Pythagorean theorem
\[
CB^2 = DC^2 + DB^2
\]

to find
\begin{equation}
a^2 = (b - c \cos \alpha)^2 + c^2 \sin^2 \alpha. \tag{24} \label{eq:24}
\end{equation}

With some algebra, we get
\begin{equation}
a^2 = b^2 + c^2\,(\cos^2 \alpha + \sin^2 \alpha) - 2bc \cos \alpha. \tag{25} \label{eq:25}
\end{equation}

Then, using the trigonometric identity $\cos^2 \alpha + \sin^2 \alpha = 1$, we obtain
\begin{equation}
a^2 = b^2 + c^2 - 2bc \cos \alpha, \tag{26} \label{eq:26}
\end{equation}
which is the Cosine Rule.

We have derived the Sine Rule [Eq.~\eqref{eq:23}] and the Cosine Rule [Eq.~\eqref{eq:26}] using the Primary Gasing Triangle and the appropriate scaling factors. In solving trigonometric problems, we recommend not relying solely on memorized rules but rather encouraging students to solve problems using the Primary Gasing Triangle and appropriate scaling factors, as explained in several examples in \hyperref[appendix1]{Appendix 1}.

\break
\section{Cofunction Identity}
\label{sec:Cofunction Identity}
In Fig.~\ref{Figure_6}, let’s name the angle $\angle EAB = \gamma = 90^\circ - \alpha$. Since $\triangle ABC$ is a Primary Gasing Triangle, $AC = \sin \gamma$ and the length of side $BC = \cos \gamma$. The lengths of the other sides can be obtained in the same manner as we did for the side lengths of the triangle in Fig.~\ref{Figure_6}, but for the angle $\gamma$. The results are as shown in Fig.~\ref{Figure_10}.

% Reset figure counter to 9 so that Figure 10 is correctly numbered
\setcounter{figure}{9}
\begin{figure}[H]
    \centering
    \scalebox{1}{\input{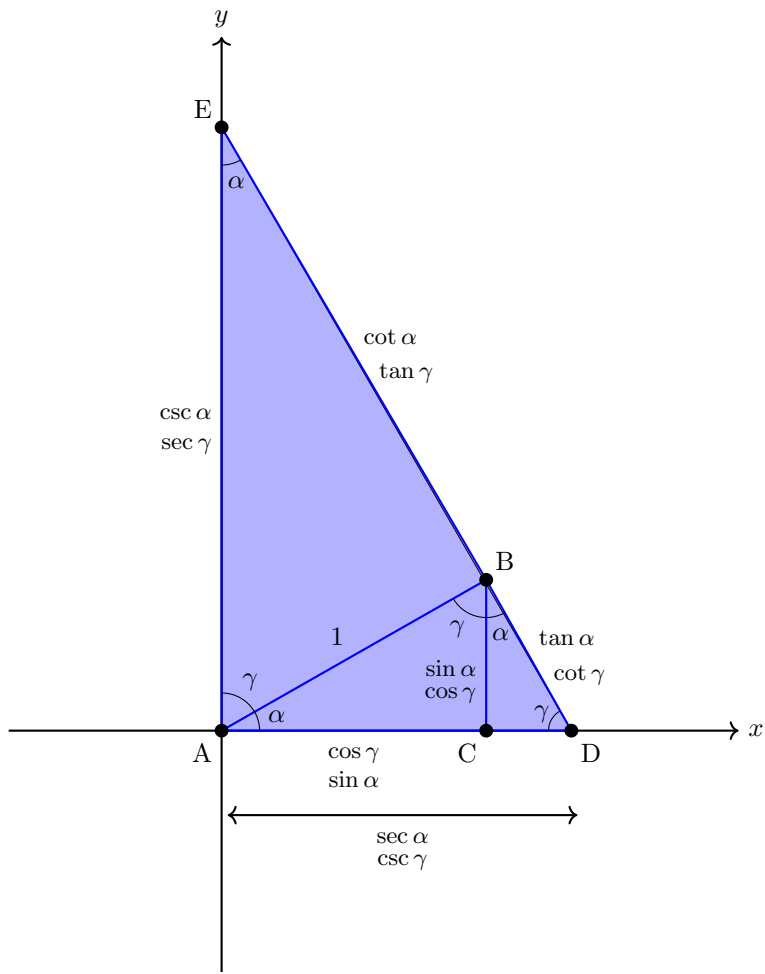}}
    \caption{Triangles to derive the Cofunction Identity}
    \label{Figure_10}
\end{figure}

From Fig.~\ref{Figure_10}, we observe the cofunction identities for the six trigonometric functions:
\begin{align}
\cos \gamma &= \cos \left(\frac{\pi}{2} - \alpha\right) = \sin \alpha \tag{27} \label{eq:27} \\
\sin \gamma &= \sin \left(\frac{\pi}{2} - \alpha\right) = \cos \alpha \tag{28} \label{eq:28} \\
\sec \gamma &= \sec \left(\frac{\pi}{2} - \alpha\right) = \csc \alpha \tag{29} \label{eq:29} \\
\csc \gamma &= \csc \left(\frac{\pi}{2} - \alpha\right) = \sec \alpha \tag{30} \label{eq:30} \\
\tan \gamma &= \tan \left(\frac{\pi}{2} - \alpha\right) = \cot \alpha \tag{31} \label{eq:31} \\
\cot \gamma &= \cot \left(\frac{\pi}{2} - \alpha\right) = \tan \alpha \tag{32} \label{eq:32}
\end{align}

These cofunction identities illustrate how each trigonometric function of an angle is related to another function of its complementary angle.

\break
\section{Primary Trigonometric Functions for Obtuse Angles}
\label{sec:Primary Trigonometric Functions for Obtuse Angles}
We have seen that the most important trigonometric functions are sine and cosine, which we refer to as the Primary Trigonometric Functions. We have also already explored the definitions of sine and cosine for acute angles as the lengths of the vertical and horizontal sides of the Primary Gasing Triangle.

For obtuse angles $(\theta > 90^\circ)$, the definition of the side lengths of the triangle becomes the projection of the unit vector onto the $y$-axis and $x$-axis in Cartesian coordinates. This definition differs from the one given in the unit circle concept ([15], [16], [17]). In the unit circle concept, the values of the trigonometric functions are determined solely based on the coordinates of the point on the terminal side of the angle that lies on the unit circle, without providing a physical meaning for those coordinates.

\begin{figure}[H]
    \centering
    \scalebox{1.1}{\input{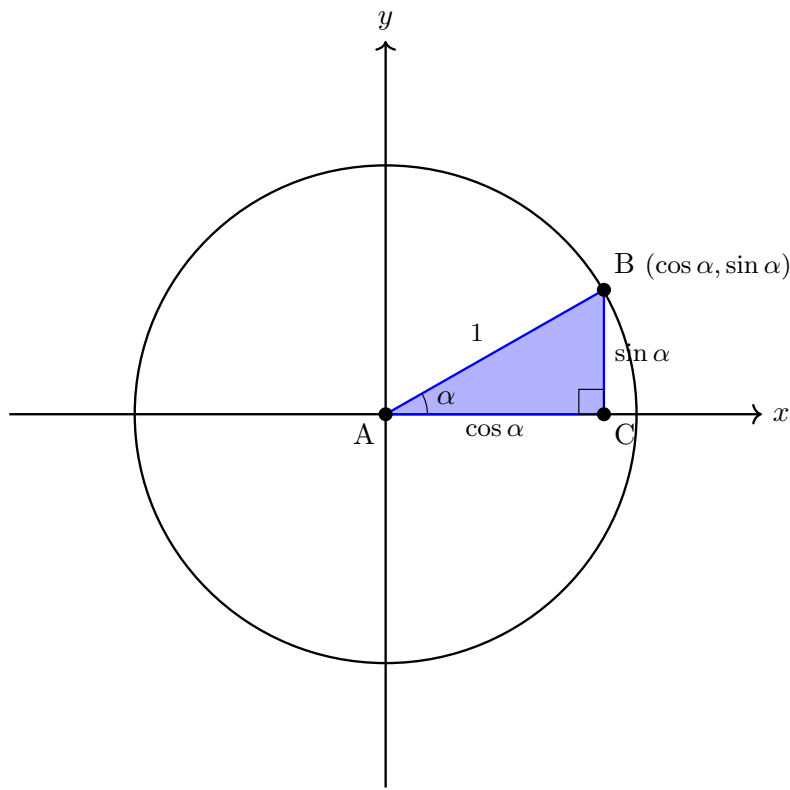}}
    \renewcommand{\thefigure}{11a} % Temporarily redefine figure numbering
    \caption{Primary Gasing Triangle in quadrant 1}
    \label{Figure_11a}
    \renewcommand{\thefigure}{\arabic{figure}} % Restore standard numbering
\end{figure}

Assume the Primary Gasing Triangle $\triangle ABC$ in Fig.~\ref{Figure_1} is placed on the Cartesian coordinate system such that side $AC$ lies on the $x$-axis and $\angle BAC = \alpha$. Here, the length of $AC$ is $\cos \alpha$ and the length of $CB$ is $\sin \alpha$. The coordinates of point $B$ are $B\,(\cos \alpha, \sin \alpha)$, as shown in Fig.~\ref{Figure_11a}.

\begin{figure}[H]
    \centering
    \scalebox{1.1}{\input{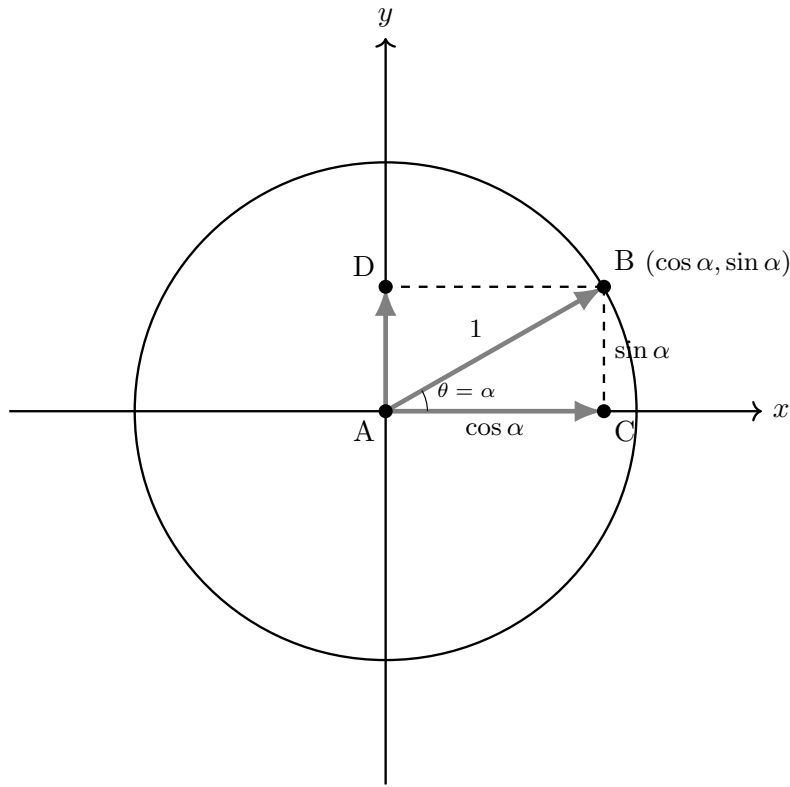}}
    \renewcommand{\thefigure}{11b} % Temporarily redefine figure numbering
    \caption{Transforming side \( AB \) of the Primary Gasing Triangle into a unit vector forming an angle \( \theta \) with the \( x \)-axis, where \( \theta = \alpha \)}
    \label{Figure_11b}
    \renewcommand{\thefigure}{\arabic{figure}} % Restore standard numbering
\end{figure}

Next, we transform side $AB$ of this triangle into a unit vector that forms an angle $\theta$ with the $x$-axis. Here, $\theta = \alpha$, as shown in Fig.~\ref{Figure_11b}.

We define $\sin\,\theta$ as the projection of the unit vector $AB$ onto the $y$-axis, while $\cos\,\theta$ is the projection of the unit vector $AB$ onto the $x$-axis. Here, the coordinates of point $B$ represent these projections: the abscissa ($x$-coordinate) represents the projection of the unit vector $AB$ onto the $x$-axis, and the ordinate ($y$-coordinate) represents the projection onto the $y$-axis. The unit vector $AB$ can rotate freely around point $A$.

Using the new definitions of sine and cosine, let's calculate the values for all obtuse angles.

\subsection{Obtuse Angle in the Second Quadrant $(90^\circ < \theta < 180^\circ)$}
\begin{enumerate}[label=\alph*), left=12pt, itemsep=1em]
    \item \textbf{Rotate the Unit Vector} \\
    Rotate unit vector $AB = 1$ unit to form an angle $\theta = (180^\circ - \alpha)$ with the $x$-axis, as shown in Fig.~\ref{Figure_11c}.

    \begin{figure}[H]
        \centering
        \scalebox{1.1}{\input{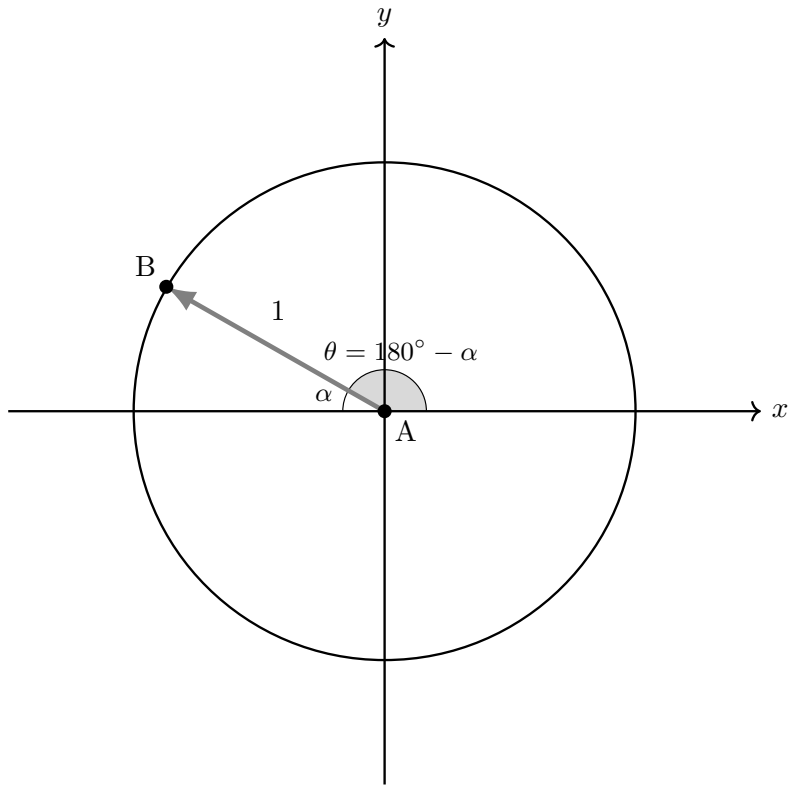}}
        \renewcommand{\thefigure}{11c} % Temporarily redefine figure numbering
        \caption{Unit vector $AB$ in quadrant 2}
        \label{Figure_11c}
        \renewcommand{\thefigure}{\arabic{figure}} % Restore standard numbering
    \end{figure}
    
    \item \textbf{Construct the Primary Gasing Triangle $\triangle ABC$} \\
    As shown in Fig.~\ref{Figure_11d}, with the Primary Gasing Triangle properties, we have $AC = \cos \alpha$ and $CB = \sin \alpha$.

    \begin{figure}[H]
        \centering
        \scalebox{1.1}{\input{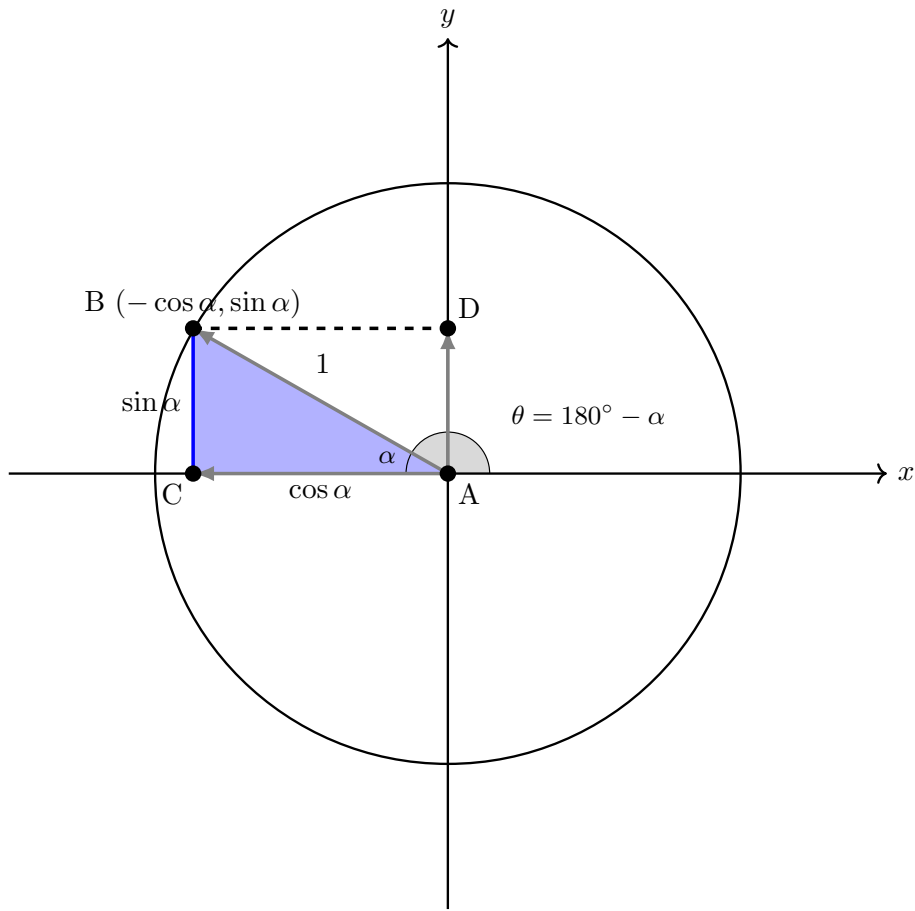}}
        \renewcommand{\thefigure}{11d} % Temporarily redefine figure numbering
        \caption{Calculating sine and cosine of an angle in quadrant 2}
        \label{Figure_11d}
        \renewcommand{\thefigure}{\arabic{figure}} % Restore standard numbering
    \end{figure}
    
    \item \textbf{Determine the Coordinates of Point $B$} \\
    The coordinates of point $B$ are $B\,(-\cos \alpha, \sin \alpha)$.

    \item \textbf{Calculate Projections} \\
    The value of $\sin\,(180^\circ - \alpha)$ is the projection of the unit vector $AB$ onto the $y$-axis, which is the vector $AD$, equal to $+\sin \alpha$ (positive because the unit vector is directed towards the positive $y$-axis). This is the same as the ordinate ($y$-coordinate) of point $B$.
    
    The value of $\cos\,(180^\circ - \alpha)$ is the projection of the unit vector $AB$ onto the $x$-axis, which is the vector $AC$, equal to $-\cos \alpha$ (negative because the unit vector is directed towards the negative $x$-axis). This is the same as the abscissa ($x$-coordinate) of point $B$.
\end{enumerate}

We can write the results as:
\begin{align}
\cos\,(180^\circ - \alpha) &= -\cos \alpha \tag{33} \label{eq:33} \\
\sin\,(180^\circ - \alpha) &= \sin \alpha \tag{34} \label{eq:34}
\end{align}

These results show how the values of sine and cosine are determined for an obtuse angle in the second quadrant using the projection method with the Gasing Triangle. This approach maintains a clear physical interpretation by using Cartesian coordinates and projections.

\subsection{Obtuse Angle in the Third Quadrant $(180^\circ < \theta < 270^\circ)$}

\begin{figure}[H]
    \centering
    \scalebox{1.1}{\input{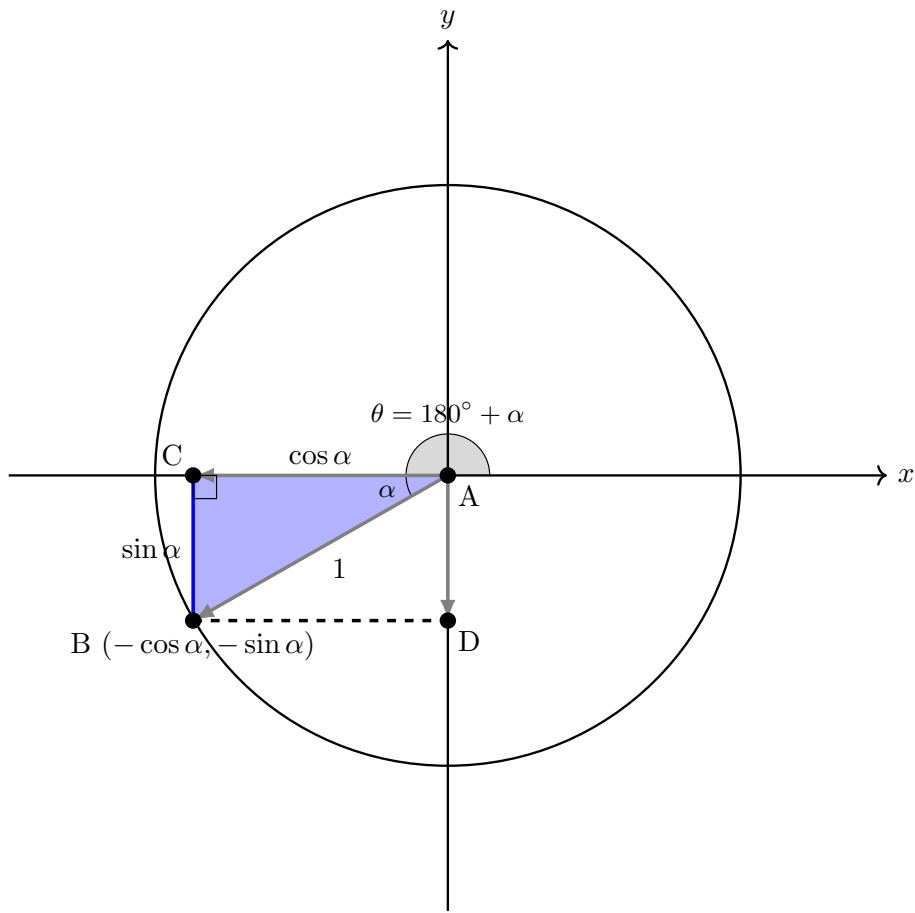}}
    \renewcommand{\thefigure}{11f} % Temporarily redefine figure numbering
    \caption{Calculating sine and cosine of an angle in quadrant 3}
    \label{Figure_11f}
    \renewcommand{\thefigure}{\arabic{figure}} % Restore standard numbering
\end{figure}

Using the steps we followed for obtuse angles in the second quadrant, and referring to Figs.~\ref{Figure_11e} and~\ref{Figure_11f}, we obtain the following results:

\begin{figure}[H]
    \centering
    \scalebox{1.1}{\input{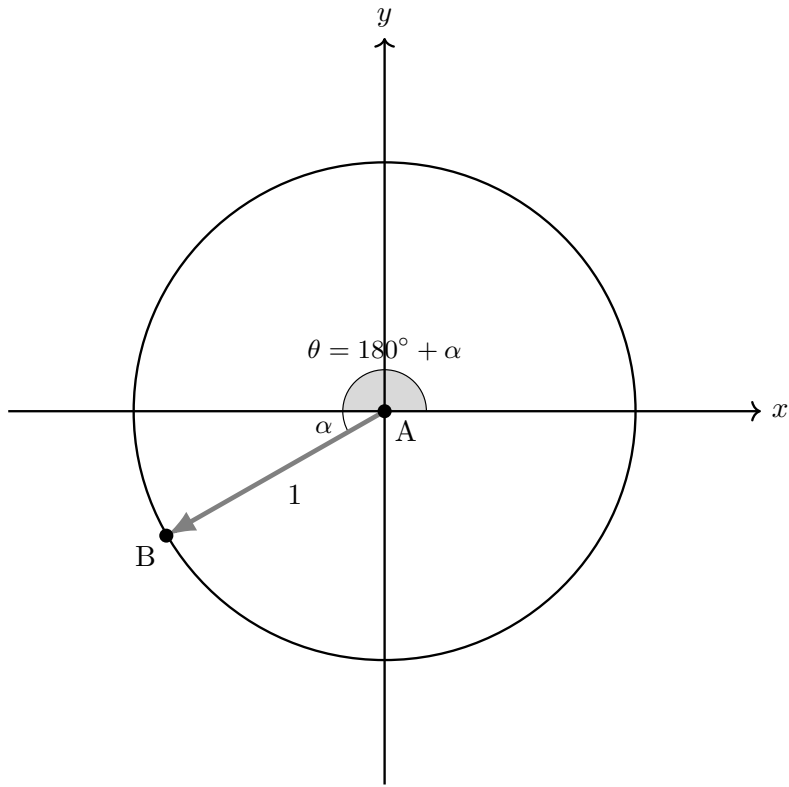}}
    \renewcommand{\thefigure}{11e} % Temporarily redefine figure numbering
    \caption{Unit vector $AB$ in quadrant 3}
    \label{Figure_11e}
    \renewcommand{\thefigure}{\arabic{figure}} % Restore standard numbering
\end{figure}

The value of $\sin\,(180^\circ + \alpha)$ is the projection of the unit vector $AB$ onto the $y$-axis, which is the vector $AD$, equal to $-\sin \alpha$ (negative because the unit vector is directed towards the negative $y$-axis). This is the same as the ordinate ($y$-coordinate) of point $B$.

The value of $\cos\,(180^\circ + \alpha)$ is the projection of the unit vector $AB$ onto the $x$-axis, which is the vector $AC$, equal to $-\cos \alpha$ (negative because the unit vector is directed towards the negative $x$-axis). This is the same as the abscissa ($x$-coordinate) of point $B$.

We can write the results as:
\begin{align}
\cos\,(180^\circ + \alpha) &= -\cos \alpha \tag{35} \label{eq:35} \\
\sin\,(180^\circ + \alpha) &= -\sin \alpha \tag{36} \label{eq:36}
\end{align}

\subsection{Obtuse Angle in the Fourth Quadrant $(270^\circ < \theta < 360^\circ)$}

\begin{figure}[H]
    \centering
    \scalebox{1.1}{\input{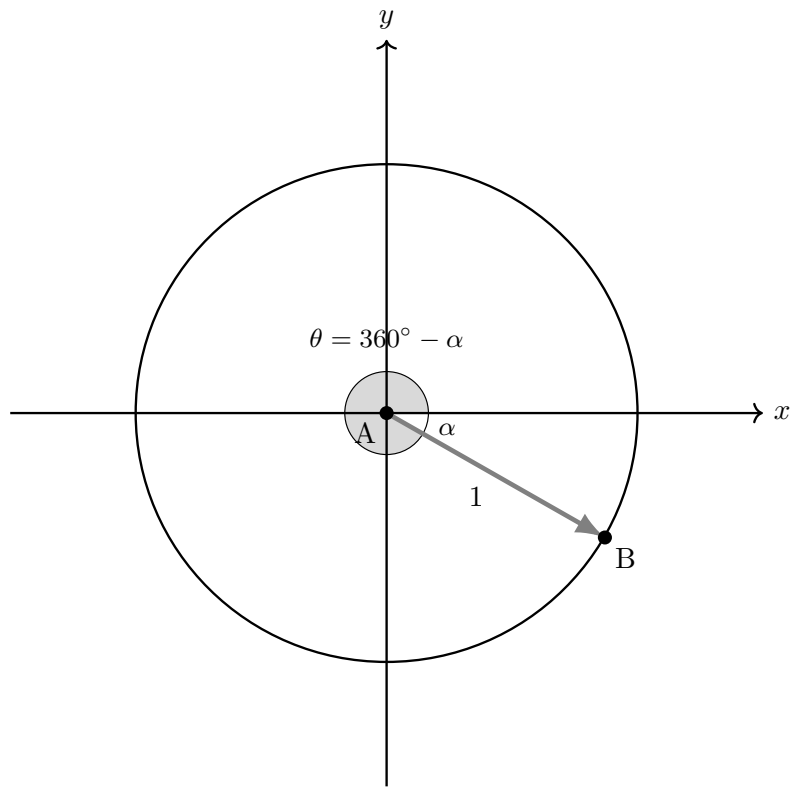}}
    \renewcommand{\thefigure}{11g} % Temporarily redefine figure numbering
    \caption{Unit vector $AB$ in quadrant 4}
    \label{Figure_11g}
    \renewcommand{\thefigure}{\arabic{figure}} % Restore standard numbering
\end{figure}

\begin{figure}[H]
    \centering
    \scalebox{1.1}{\input{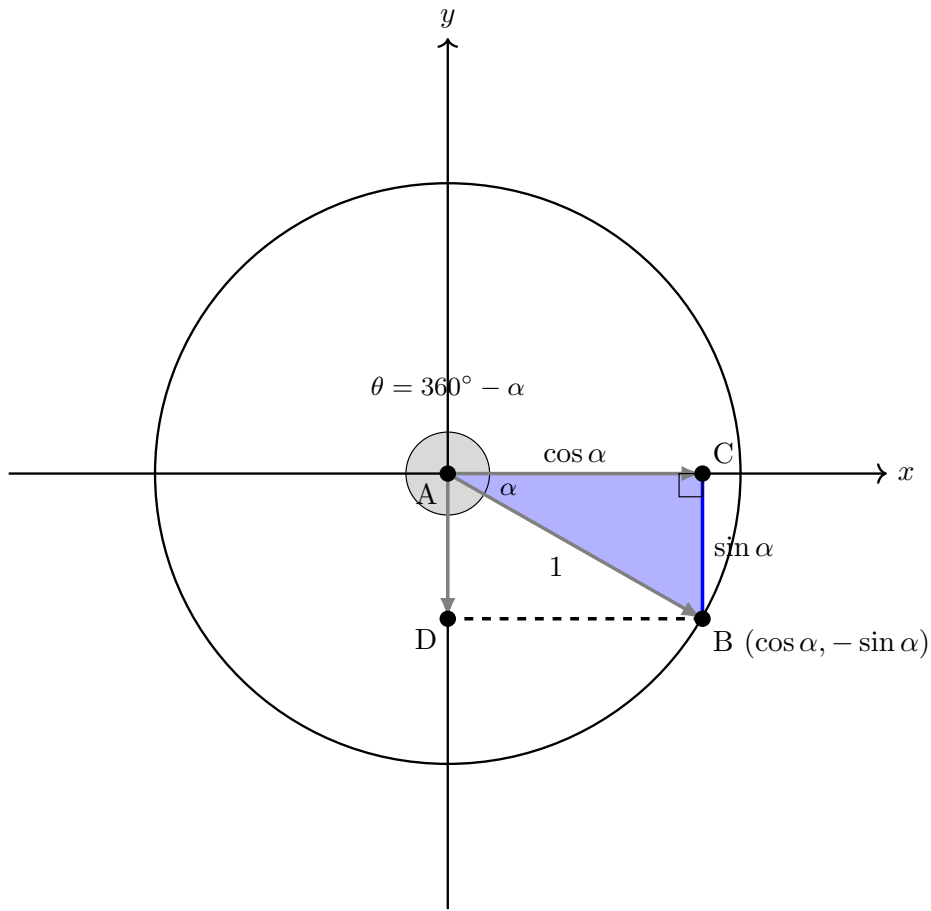}}
    \renewcommand{\thefigure}{11h} % Temporarily redefine figure numbering
    \caption{Calculating sine and cosine of an angle in quadrant 4}
    \label{Figure_11h}
    \renewcommand{\thefigure}{\arabic{figure}} % Restore standard numbering
\end{figure}

Using the steps we followed for obtuse angles in the second quadrant, and referring to Figs.~\ref{Figure_11g} and~\ref{Figure_11h}, we obtain the following results:

\begin{align}
\cos\,(360^\circ - \alpha) &= \cos \alpha \tag{37} \label{eq:37} \\
\sin\,(360^\circ - \alpha) &= -\sin \alpha \tag{38} \label{eq:38}
\end{align}

Defining sine and cosine as the projections of a unit vector in Cartesian coordinates allows us to calculate the values of these trigonometric functions for any angle of inclination of this unit vector. This enables these trigonometric functions to describe periodic motion, such as vibrations, waves, light, etc. 

By using this method, we can consistently extend the definitions of sine and cosine to cover all angles, maintaining a clear physical interpretation of their values based on projections in the Cartesian plane.

\break
\section{Conclusion}
\label{sec:conclusion}
The primary trigonometric functions $\sin \theta$ and $\cos \theta$ are defined as the projections on the $y$-axis and $x$-axis of a unit vector that forms an angle $\theta$ with the $x$-axis. 

For acute angles, the projections of the unit vector, $\sin \theta$ and $\cos\,\theta$, can be reduced to the side lengths of the Primary Gasing Triangle, as shown in Fig.~\ref{Figure_2}.

For other trigonometric functions, or derived trigonometric functions (secant, tangent, cosecant, and cotangent), their values are determined based on the primary trigonometric functions, given by Eqs.~\eqref{eq:1}, \eqref{eq:2}, \eqref{eq:3}, and \eqref{eq:4}. For acute angles, these derived trigonometric functions are represented as the side lengths of the derived Gasing Triangles shown in Figs.~\ref{Figure_5b} and~\ref{Figure_5c}.

All derivations of formulas and solutions to trigonometric problems can be accomplished using the Primary Gasing Triangle and the appropriate scaling factors. This makes learning trigonometry simpler, more understandable, and more intuitive. Students do not need to memorize numerous formulas that may not be particularly useful later.

With the Primary Gasing Triangle, we have also demonstrated that it is quite feasible to prove the Pythagorean Theorem using trigonometry.

The definition of Primary Trigonometric Functions for all angles allows trigonometric functions to better explain various periodic phenomena such as vibration motion, waves, sound, light, or other aspects related to periodic motion. This will make learning trigonometry even more meaningful, and students will be more likely to enjoy the subject.

\break
\appendix
\renewcommand{\thesection}{}  % Remove automatic numbering (A, B, C, etc.)
\section*{Appendix 1: Solving Trigonometry Problems with the Gasing Triangle}
\label{appendix1}

\subsection*{Problem 1:}
If \(\sin \alpha = \frac{1}{2}\), what is \(\cos \alpha\)? \newline
\break
\noindent \textbf{Solution:}

\begin{figure}[H]
    \centering
    \scalebox{1.4}{\input{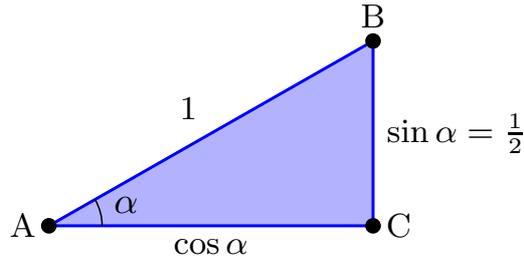}}
    \renewcommand{\thefigure}{A1.1} % Temporarily redefine figure numbering
    \captionsetup{labelformat=empty} % Removes default label format
    \caption{Figure A1.1: Primary Gasing Triangle with \( \sin \alpha = \frac{1}{2} \)}
    \label{Figure_A1.1}
    \renewcommand{\thefigure}{\arabic{figure}} % Restore standard numbering
\end{figure}

Let's draw the Primary Gasing Triangle as shown in Fig.~\ref{Figure_A1.1}.

We can calculate \(\cos \alpha\) using the Pythagorean Theorem in \(\triangle ABC\):
\[
\cos \alpha = \sqrt{1^2 - \left(\frac{1}{2}\right)^2} = \frac{1}{2}\sqrt{3}
\]

This problem can also be solved using the trigonometric identity \(1 = \cos^2 \alpha + \sin^2 \alpha\). However, solving it by memorizing formulas does not train students to think critically and logically, which is the spirit of teaching. By using the Primary Gasing Triangle, students can visually and logically understand the relationships between the sides and angles, leading to a deeper comprehension of trigonometry.

\subsection*{Problem 2:}
If \(\tan \alpha = \frac{3}{4}\), what is \(\cos \alpha\)? \newline
\break
\textbf{Solution:}

\begin{figure}[H]
    \centering
    \scalebox{1.4}{\input{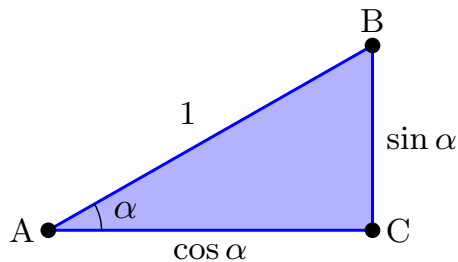}}
    \renewcommand{\thefigure}{A1.2a} % Temporarily redefine figure numbering
    \captionsetup{labelformat=empty} % Removes default label format
    \caption{Figure A1.2a: Primary Gasing Triangle with angle \(\alpha\)}
    \label{Figure_A1.2a}
    \renewcommand{\thefigure}{\arabic{figure}} % Restore standard numbering
\end{figure}

First, let's create the Primary Gasing Triangle \(\triangle ABC\) with \(\angle BAC = \alpha\) (Fig.~\ref{Figure_A1.2a}).

\begin{figure}[H]
    \centering
    \scalebox{1.4}{\input{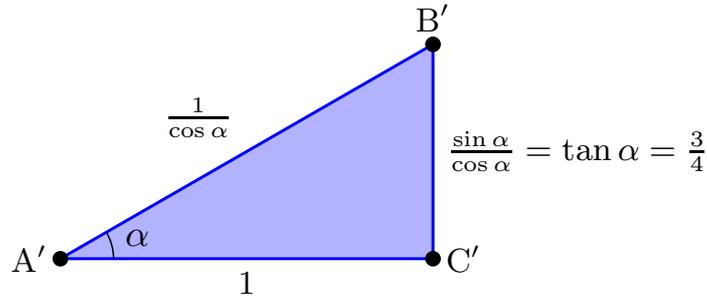}}
    \renewcommand{\thefigure}{A1.2b} % Temporarily redefine figure numbering
    \captionsetup{labelformat=empty} % Removes default label format
    \caption{Figure A1.2b: Calculating the side lengths of Triangle \( A'B'C' \) using the scale factor $\frac{1}{\cos \alpha}$ on the triangle in Figure A1.2a}
    \label{Figure_A1.2b}
    \renewcommand{\thefigure}{\arabic{figure}} % Restore standard numbering
\end{figure}

Then, use the scaling factor \(\frac{1}{\cos \alpha}\) to obtain Fig.~\ref{Figure_A1.2b}.

Using the Pythagorean Theorem in \(\triangle A'B'C'\), we get
\[
\frac{1}{\cos \alpha} = \sqrt{1^2 + \left(\frac{3}{4}\right)^2} = \frac{5}{4}. 
\]

Therefore,
\[
\cos \alpha = \frac{4}{5}.
\]

Again, by using the Primary Gasing Triangle approach, we can visually and logically determine the relationships between the trigonometric functions, leading to a deeper understanding of trigonometry concepts.

\subsection*{Problem 3:}
Calculate the value of \(a\) in Fig.~\ref{Figure_A1.3a} below: 

\begin{figure}[H]
    \centering
    \scalebox{0.5}{\input{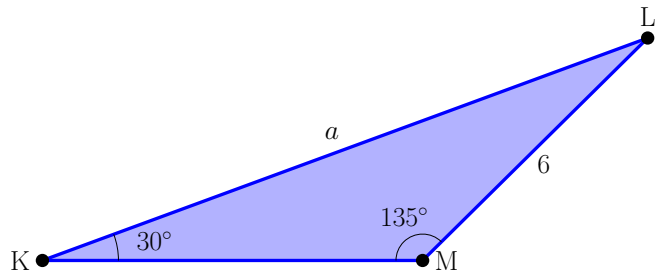}}
    \renewcommand{\thefigure}{A1.3a} % Temporarily redefine figure numbering
    \captionsetup{labelformat=empty} % Removes default label format
    \caption{Figure A1.3a: Determining the value of \( a \) in Triangle \( KLM \)}
    \label{Figure_A1.3a}
    \renewcommand{\thefigure}{\arabic{figure}} % Restore standard numbering
\end{figure}

\noindent \textbf{Solution:}

The first step is to create the Primary Gasing Triangles \(\triangle ABC\) (Fig.~\ref{Figure_A1.3b}) and \(\triangle DEF\) (Fig.~\ref{Figure_A1.3c}). Then use the scaling factor $a$ on the Primary Gasing Triangle \(\triangle ABC\) to obtain triangle \(\triangle KLN\). 

\begin{figure}[H]
    \centering
    \scalebox{0.5}{\input{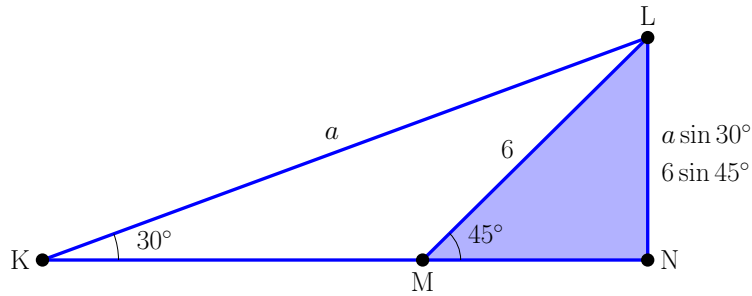}}
    \renewcommand{\thefigure}{A1.3d} % Temporarily redefine figure numbering
    \captionsetup{labelformat=empty} % Removes default label format
    \caption{Figure A1.3d: Primary Gasing Triangles to calculate the value of $a$}
    \label{Figure_A1.3d}
    \renewcommand{\thefigure}{\arabic{figure}} % Restore standard numbering
\end{figure}

Similarly, use the scaling factor $6$ on the Primary Gasing Triangle $\triangle DEF$ to obtain triangle $\triangle MLN$, as shown in Fig.~\ref{Figure_A1.3d}.

\begin{figure}[H]
    \centering
    \begin{minipage}{0.45\textwidth}
        \centering
        \scalebox{1.4}{\input{figures/Figure_A1.3b}}
        \renewcommand{\thefigure}{A1.3b} % Temporarily redefine figure numbering
        \captionsetup{labelformat=empty} % Removes default label format
        \caption{Figure A1.3b: Primary Gasing Triangle with angle $30^\circ$}
        \label{Figure_A1.3b}
        \renewcommand{\thefigure}{\arabic{figure}} % Restore standard numbering
    \end{minipage}
    \hfill
    \begin{minipage}{0.45\textwidth}
        \centering
        \scalebox{1.4}{\input{figures/Figure_A1.3c}}
        \renewcommand{\thefigure}{A1.3c} % Temporarily redefine figure numbering
        \captionsetup{labelformat=empty} % Removes default label format
        \caption{Figure A1.3c: Primary Gasing Triangle with angle $45^\circ$}
        \label{Figure_A1.3c}
        \renewcommand{\thefigure}{\arabic{figure}} % Restore standard numbering
    \end{minipage}
\end{figure}

On side $LN$, we observe that $a \sin 30^\circ = 6 \sin 45^\circ$.  
Therefore,
\[
a \cdot \frac{1}{2} = 6 \cdot \frac{1}{2} \sqrt{2}
\]
which results in
\[
a = 6 \sqrt{2}.
\]

Why don't we just use the Sine Rule from Eq.~\eqref{eq:23}? It's straightforward. That's true, but again we are not too interested in having students memorize formulas. Instead, we want to encourage students to prioritize creativity, computational logic, and critical thinking. This approach makes learning trigonometry more engaging.

\subsection*{Problem 4:}
\begin{figure}[H]
    \centering
    \scalebox{0.4}{\input{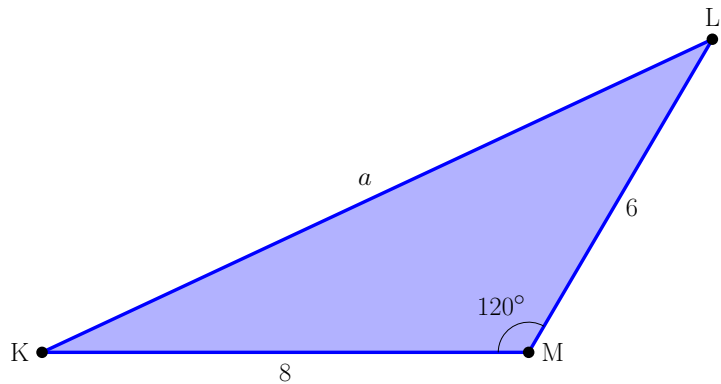}}
    \renewcommand{\thefigure}{A1.4a} % Temporarily redefine figure numbering
    \captionsetup{labelformat=empty} % Removes default label format
    \caption{Fig.~A1.4a: Determining the value of \( a \) in Triangle \( KLM \) with \( \angle KML = 120^\circ \)}
    \label{Figure_A1.4a}
    \renewcommand{\thefigure}{\arabic{figure}} % Restore standard numbering
\end{figure}

Calculate the value of \(a\) in the triangle in Fig.~\ref{Figure_A1.4a}.

\noindent \textbf{Solution:}

We can solve this using the Cosine Rule as shown in Eq.~\eqref{eq:26}. However, in the spirit of Gasing, we prefer not to have students memorize formulas. A more engaging method is to use the Primary Gasing Triangle and the associated scaling factor. \newline

\begin{figure}[H]
    \centering
    \scalebox{1.2}{\input{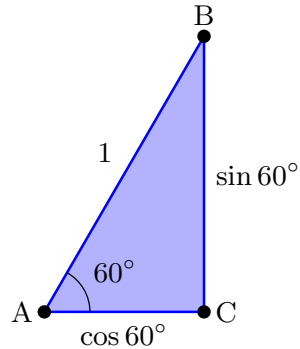}}
    \renewcommand{\thefigure}{A1.4b} % Temporarily redefine figure numbering
    \captionsetup{labelformat=empty} % Removes default label format
    \caption{Fig.~A1.4b: Primary Gasing Triangle with angle $60^\circ$}
    \label{Figure_A1.4b}
    \renewcommand{\thefigure}{\arabic{figure}} % Restore standard numbering
\end{figure}

\begin{figure}[H]
    \centering
    \scalebox{0.4}{\input{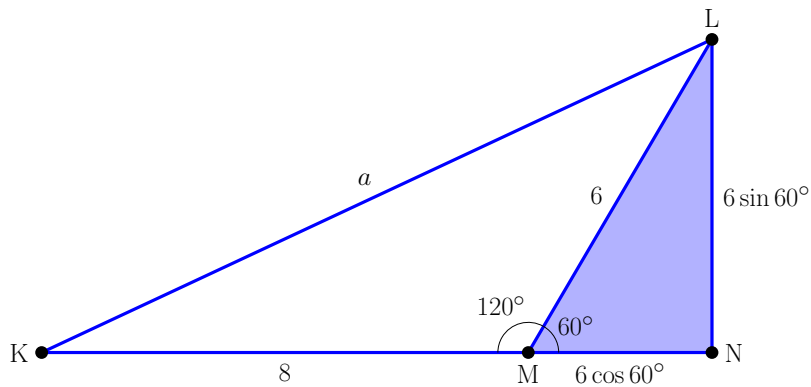}}
    \renewcommand{\thefigure}{A1.4c} % Temporarily redefine figure numbering
    \captionsetup{labelformat=empty} % Removes default label format
    \caption{Fig.~A1.4c: Primary Gasing Triangles to calculate the value of $a$}
    \label{Figure_A1.4c}
    \renewcommand{\thefigure}{\arabic{figure}} % Restore standard numbering
\end{figure}

We create the Primary Gasing Triangle \(\triangle ABC\) with angle \(\angle BAC = 60^\circ\)(Fig.~\ref{Figure_A1.4b}). Then, use the scaling factor 6 to obtain triangle \(\triangle MLN\), as shown in Fig.~\ref{Figure_A1.4c}.

Using the Pythagorean Theorem on \(\triangle KLN\):
\[
a = \sqrt{(8 + 6 \cos\,60^\circ)^2 + (6 \sin\,60^\circ)^2}
\]
\[
= \sqrt{(8 + 3)^2 + (6 \cdot \frac{1}{2} \sqrt{3})^2}
\]
From this, we get
\[
a = 2\sqrt{37}.
\]
Easier, isn't it?

\subsection*{Problem 5:}
From a point \(A\) above the ground, a person sees a flagpole \(DB\) that is 12 meters high with angles of view of \(45^\circ\) and \(30^\circ\). The flagpole is located on a hill. Calculate the height of the hill \(CD\).

\begin{figure}[H]
    \centering
    \scalebox{0.6}{\input{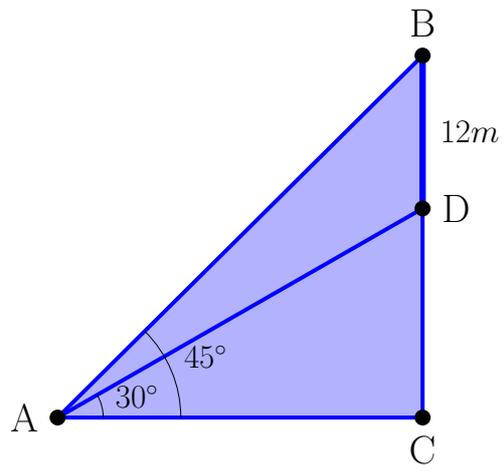}}
    \renewcommand{\thefigure}{A1.5a} % Temporarily redefine figure numbering
    \captionsetup{labelformat=empty} % Removes default label format
    \caption{Fig.~A1.5a: Determining the length of $CD$}
    \label{Figure_A1.5a}
    \renewcommand{\thefigure}{\arabic{figure}} % Restore standard numbering
\end{figure}

\break
\noindent \textbf{Solution:}

We first create the Primary Gasing Triangles \(\triangle AEF\) with angle \(\alpha = 45^\circ\) (Fig.~\ref{Figure_A1.5b}) and \(\triangle AGH\) with angle \(\alpha = 30^\circ\) (Fig.~\ref{Figure_A1.5c}).

\begin{figure}[H]
    \centering
    \begin{minipage}{0.45\textwidth}
        \centering
        \scalebox{1.4}{\input{figures/Figure_A1.5b}}
        \renewcommand{\thefigure}{A1.5b} % Temporarily redefine figure numbering
        \captionsetup{labelformat=empty} % Removes default label format
        \caption{Fig.~A1.5b: Primary Gasing Triangle with angle $45^\circ$}
        \label{Figure_A1.5b}
        \renewcommand{\thefigure}{\arabic{figure}} % Restore standard numbering
    \end{minipage}
    \hfill
    \begin{minipage}{0.45\textwidth}
        \centering
        \scalebox{1.4}{\input{figures/Figure_A1.5c}}
        \renewcommand{\thefigure}{A1.5c} % Temporarily redefine figure numbering
        \captionsetup{labelformat=empty} % Removes default label format
        \caption{Fig.~A1.5c: Primary Gasing Triangle with angle $30^\circ$}
        \label{Figure_A1.5c}
        \renewcommand{\thefigure}{\arabic{figure}} % Restore standard numbering
    \end{minipage}
\end{figure}

In Fig.~\ref{Figure_A1.5a}, triangles \(\triangle ABC\) and \(\triangle ADC\) have the same horizontal side \(AC\). Therefore, we will scale both Primary Gasing Triangles \(\triangle AEF\) and \(\triangle AGH\) so that the lengths of their horizontal sides will be equal. For this, the scaling factor for \(\triangle AEF\) is \(\frac{x}{\cos\,45^\circ}\), and the scaling factor for \(\triangle AGH\) is \(\frac{x}{\cos\,30^\circ}\).

Thus, we obtain the following figure:
\begin{figure}[H]
    \centering
    \scalebox{0.7}{\input{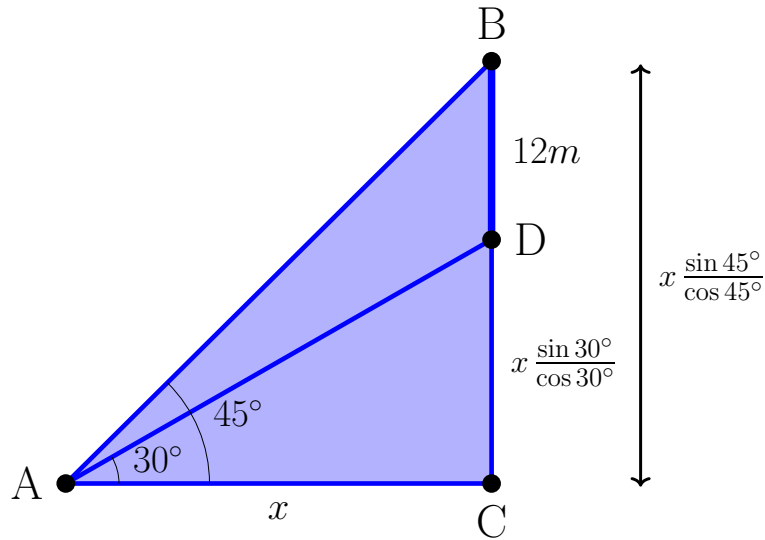}}
    \renewcommand{\thefigure}{A1.5d} % Temporarily redefine figure numbering
    \captionsetup{labelformat=empty} % Removes default label format
    \caption{Fig.~A1.5d: Calculating the value of $x$}
    \label{Figure_A1.5d}
    \renewcommand{\thefigure}{\arabic{figure}} % Restore standard numbering
\end{figure}

From the figure, we obtain:
\[
CD = x \frac{\sin 30^\circ}{\cos 30^\circ} = \frac{x}{\sqrt{3}}
\]
\[
CB = x \frac{\sin 45^\circ}{\cos 45^\circ} = x
\]

Next, since the length of \(DB\) is 12 meters, then
\[
DB = CB - CD
\]
\[
12 = x - \frac{x}{\sqrt{3}}
\]
or
\[
x = \frac{12}{\left(1 - \frac{1}{\sqrt{3}}\right)}.
\]

Thus, the height of the hill \(CD\) is
\[
CD = \frac{x}{\sqrt{3}} = \frac{12}{\left(1 - \frac{1}{\sqrt{3}}\right) \sqrt{3}} = \frac{12}{(\sqrt{3} - 1)} \text{ meters.}
\]

\break
\renewcommand{\thesection}{}  % Remove automatic numbering (A, B, C, etc.)
\section*{Appendix 2: Some Proofs of the Pythagorean Theorem with Trigonometry}
\label{appendix2}

In this section, we will prove the Pythagorean theorem using the Primary Gasing Triangle for several cases found in various references \autocite{loomis_1968, sullivan_2016, bylen_ziegler_barnett_2011, stitz_zeager_2012, sundstrom_schlicker_2019}. We have observed that various algebraic or geometric proofs of Pythagoras can be transformed into trigonometric proofs using the Primary Gasing Triangle.

In each case, we will prove \(1 = \cos^2 \alpha + \sin^2 \alpha\). We will not use the Pythagorean theorem for this proof. This proof is equivalent to proving the Pythagorean theorem as shown in Fig.~\ref{Figure_5a}. Thus, we simultaneously prove both the trigonometric identity and the Pythagorean theorem. There is no circularity of argument.

\subsection*{Case 1:}
In this figure, we construct a square \(\square ABCD\) and within it, a smaller square \(\square EFGH\) with side length 1 unit. Here, we have four similar Primary Gasing Triangles: \(\triangle HEA\), \(\triangle GHD\), \(\triangle FGC\), and \(\triangle EFB\).

\begin{figure}[H]
    \centering
    \scalebox{1.3}{\input{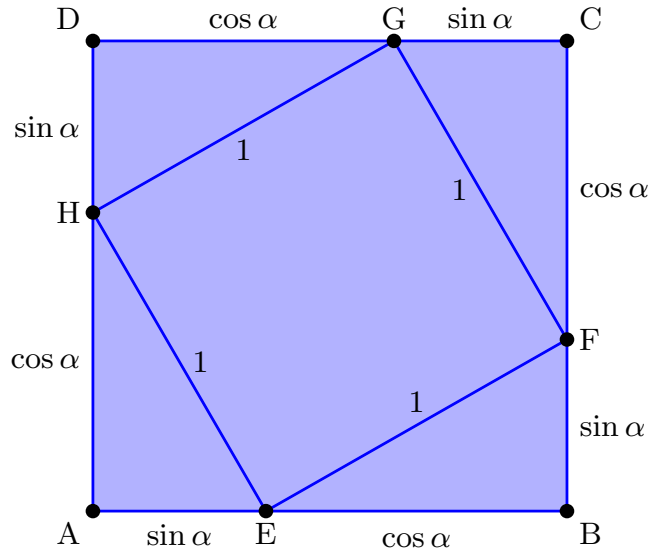}}
    \renewcommand{\thefigure}{A2.1} % Temporarily redefine figure numbering
    \captionsetup{labelformat=empty} % Removes default label format
    \caption{Fig.~A2.1: Proving the Pythagorean Theorem using the Primary Gasing Triangles \( EFB \), \( FGC \), \( GHD \), and \( HEA \)}
    \label{Figure_A2.1}
    \renewcommand{\thefigure}{\arabic{figure}} % Restore standard numbering
\end{figure}

In Fig.~\ref{Figure_A2.1}, since
\[\text{Area of } \triangle HEA = \text{Area of } \triangle GHD = \text{Area of } \triangle FGC = \text{Area of } \triangle EFB\]
and
\[\text{Area of } \square ABCD = \text{Area of } \square EFGH + 4 \times \text{Area of } \triangle EFB,\]
we can write
\[
(\cos \alpha + \sin \alpha)^2 = 1 + 4 \frac{\sin \alpha \cos \alpha}{2}
\]\[
\cos^2 \alpha + \sin^2 \alpha + 2 \sin \alpha \cos \alpha = 1 + 2 \sin \alpha \cos \alpha
\]
which reduces to
\begin{equation}
\cos^2 \alpha + \sin^2 \alpha = 1. \label{eq:A3}
\end{equation}
As we expected.

\subsection*{Case 2:}
\begin{figure}[H]
    \centering
    \scalebox{1}{\input{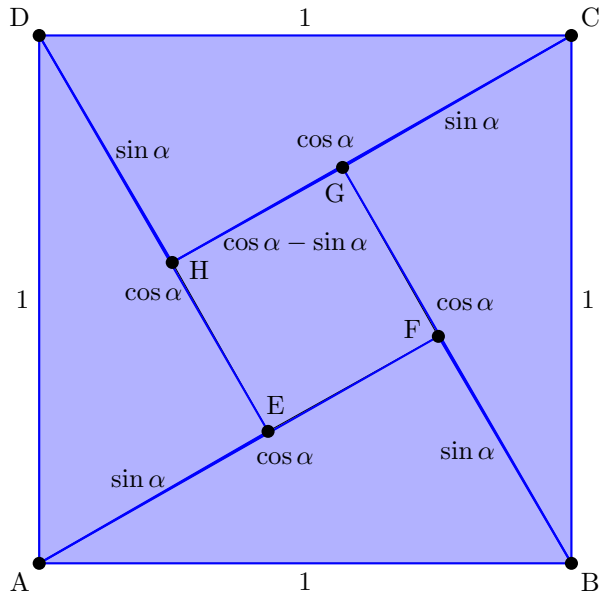}}
    \renewcommand{\thefigure}{A2.2} % Temporarily redefine figure numbering
    \captionsetup{labelformat=empty} % Removes default label format
    \caption{Fig.~A2.2: Proving the Pythagorean Theorem using a Primary Gasing Triangle different from the one in Figure A2.1}
    \label{Figure_A2.2}
    \renewcommand{\thefigure}{\arabic{figure}} % Restore standard numbering
\end{figure}

In Fig.~\ref{Figure_A2.2}, inside the square \(\square ABCD\), there is a smaller square \(\square EFGH\), and four Primary Gasing Triangles.

In the figure, it is shown that the length of \(HG\) is given by
\[
HG = \cos \alpha - \sin \alpha.
\]

Also, we observe that
\[\text{Area of } \square ABCD = \text{Area of } \square EFGH + 4 \times \text{Area of } \triangle AFB.\]
Therefore,
\[
1 = (\cos \alpha - \sin \alpha)^2 + \frac{(4 \times \sin \alpha \cos \alpha)}{2}.
\]
With a bit of algebra, we get
\[
1 = \cos^2 \alpha - 2 \cos \alpha \sin \alpha + \sin^2 \alpha + 2 \sin \alpha \cos \alpha,
\]
which simplifies to
\[
1 = \cos^2 \alpha + \sin^2 \alpha
\]
as expected.

\subsection*{Case 3:}
\begin{figure}[H]
    \centering
    \scalebox{1.1}{\input{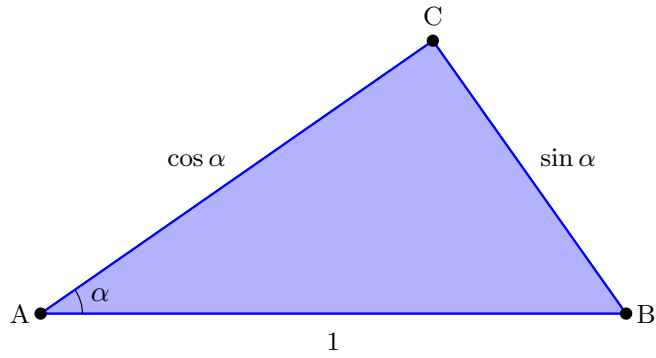}}
    \renewcommand{\thefigure}{A2.3a} % Temporarily redefine figure numbering
    \captionsetup{labelformat=empty} % Removes default label format
    \caption{Fig.~A2.3a: Primary Gasing Triangle with angle $\alpha$}
    \label{Figure_A2.3a}
    \renewcommand{\thefigure}{\arabic{figure}} % Restore standard numbering
\end{figure}

In Fig.~\ref{Figure_A2.3a}, we have a Primary Gasing Triangle \(\triangle ABC\) with the length of \(AB\) equal to 1 unit.

\begin{figure}[H]
    \centering
    \scalebox{1.1}{\input{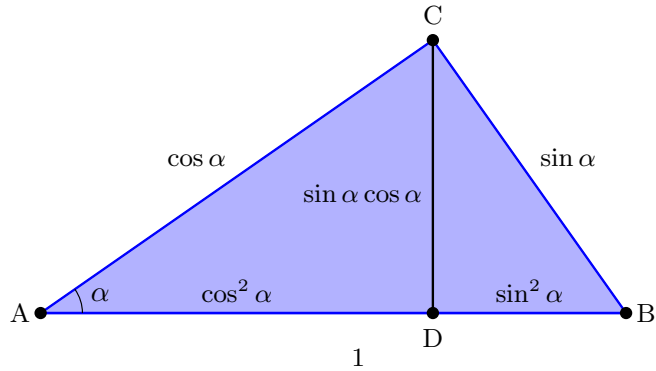}}
    \renewcommand{\thefigure}{A2.3b} % Temporarily redefine figure numbering
    \captionsetup{labelformat=empty} % Removes default label format
    \caption{Fig.~A2.3b: Proving the Pythagorean Theorem by dividing the Primary Gasing Triangle \( ACB \) into Triangles \( ACD \) and \( BCD \)}
    \label{Figure_A2.3b}
    \renewcommand{\thefigure}{\arabic{figure}} % Restore standard numbering
\end{figure}

Then, we draw a line \(CD\) perpendicular to line \(AB\). Next, use the scaling factor \(\cos \alpha\) to obtain the sides of triangle \(\triangle ACD\), where \(AD = \cos^2 \alpha\) and \(CD = \sin \alpha \cos \alpha\). Then use the scaling factor \(\sin \alpha\) to obtain \(DB = \sin^2 \alpha\), as shown in Fig.~\ref{Figure_A2.3b}.

From the figure, it can be seen that
\[
AD + DB = AB.
\]
Thus, we obtain
\[
\cos^2 \alpha + \sin^2 \alpha = 1
\]
and have proven it as expected.

\subsection*{Case 4:}
\begin{figure}[H]
    \centering
    \scalebox{1.1}{\input{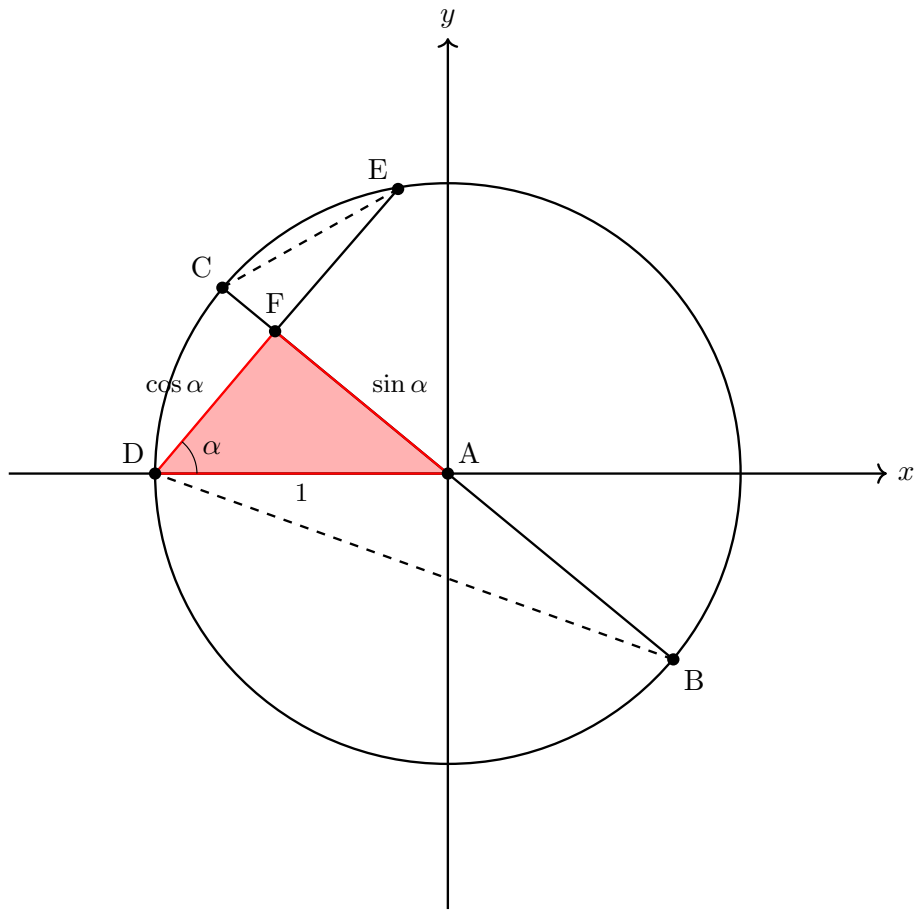}}
    \renewcommand{\thefigure}{A2.4a} % Temporarily redefine figure numbering
    \captionsetup{labelformat=empty} % Removes default label format
    \caption{Fig.~A2.4a: Proving the Pythagorean Theorem using the properties of a triangle inscribed in a circle}
    \label{Figure_A2.4a}
    \renewcommand{\thefigure}{\arabic{figure}} % Restore standard numbering
\end{figure}

Fig.~\ref{Figure_A2.4a} shows the diameter \(BC\) of a circle perpendicular to line \(DE\) where \(AF <\) radius of the circle \(AC\).

We create the Gasing Triangle \(\triangle ADF\) with hypotenuse \(AD = 1\) unit. 

Since points \(E\) and \(B\) lie on the circle, triangle \(\triangle BDF\) is similar to triangle \(\triangle ECF\). 

Since \(DE\) is perpendicular to \(CB\), the length of \(FE = DF = \cos \alpha\). 

Let \(\angle FEC = \beta\).

\begin{figure}[H]
    \centering
    \scalebox{1.4}{\input{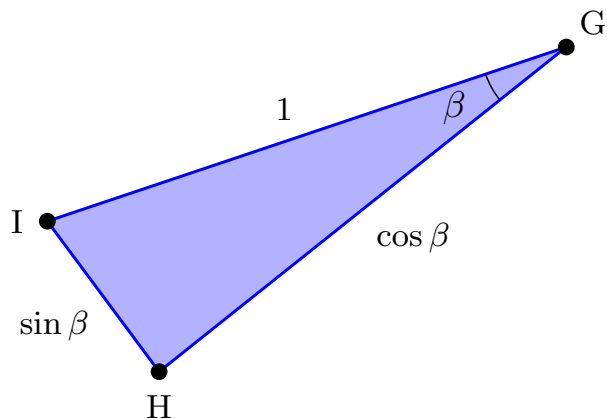}}
    \renewcommand{\thefigure}{A2.4b} % Temporarily redefine figure numbering
    \captionsetup{labelformat=empty} % Removes default label format
    \caption{Fig.~A2.4b: Primary Gasing Triangle with angle $\beta$}
    \label{Figure_A2.4b}
    \renewcommand{\thefigure}{\arabic{figure}} % Restore standard numbering
\end{figure}

Next, we create the Gasing Triangle \(\triangle GIH\) similar to triangle \(\triangle ECF\). By using the scale factor \(\frac{\cos \alpha}{\cos \beta}\) on triangle \(\triangle GIH\), we obtain triangle \(\triangle ECF\). Thus, we find the side lengths \(FC = \sin \beta \times \frac{\cos \alpha}{\cos \beta}\) and \(EC = \frac{\cos \alpha}{\cos \beta}\).

\begin{figure}[H]
    \centering
    \scalebox{0.9}{\input{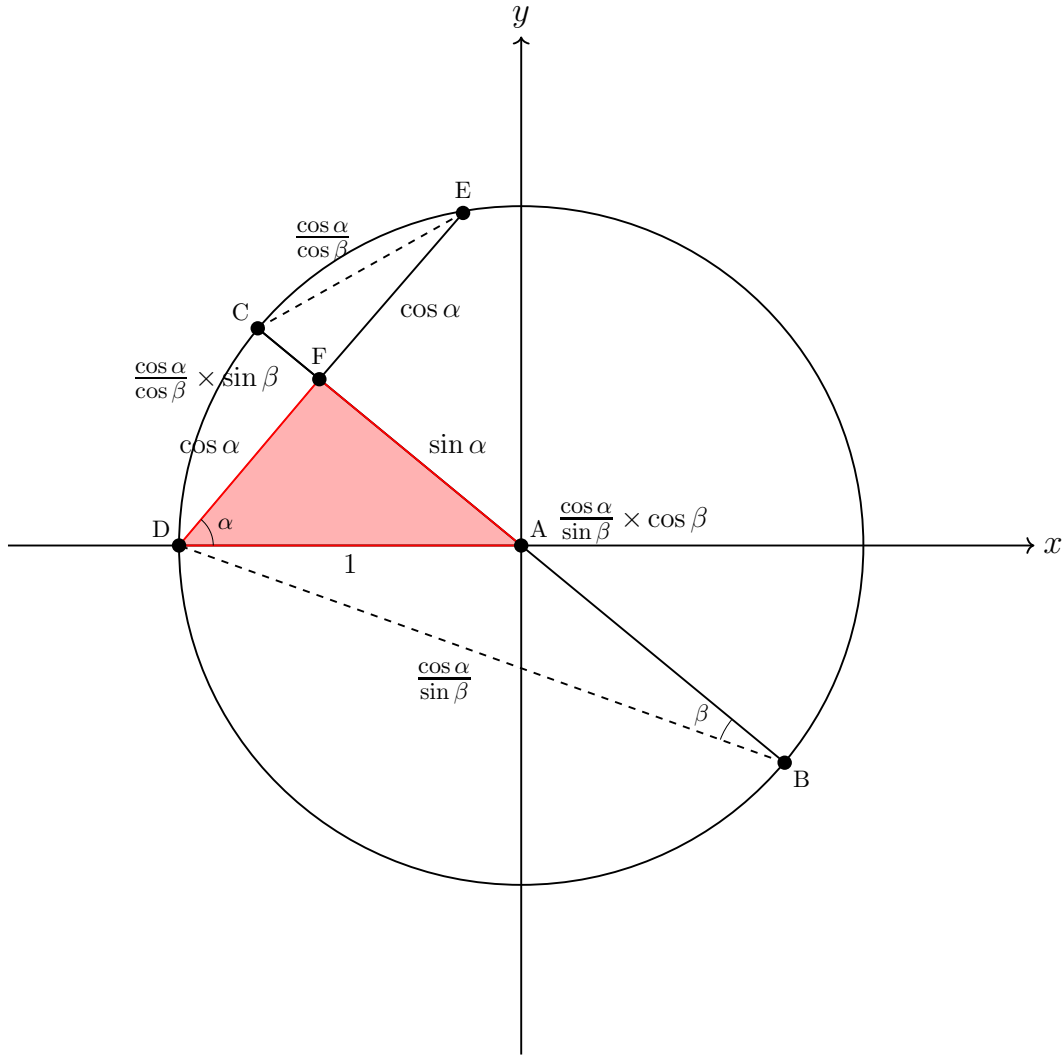}}
    \renewcommand{\thefigure}{A2.4c} % Temporarily redefine figure numbering
    \captionsetup{labelformat=empty} % Removes default label format
    \caption{Fig.~A2.4c: Side lengths of Triangles \( ECF \), \( DAF \), and \( DBF \) for proving the Pythagorean Theorem}
    \label{Figure_A2.4c}
    \renewcommand{\thefigure}{\arabic{figure}} % Restore standard numbering
\end{figure}

Next, use the scaling factor \(\frac{\cos \alpha}{\sin \beta}\) on the Gasing Triangle \(\triangle GIH\) to obtain the sides of \(\triangle BDF\), which are
\[
FB = \cos \beta \times \frac{\cos \alpha}{\sin \beta}
\]
and
\[
DB = \frac{\cos \alpha}{\sin \beta}.
\]
The complete side lengths of these triangles are shown in Fig.~\ref{Figure_A2.4c}.

Since the radius of the circle is 1 unit, then \(FC = 1 - \sin \alpha\) and \(FB = 1 + \sin \alpha\). Thus, we have the equations:
\[
1 - \sin \alpha = \frac{\cos \alpha}{\cos \beta} \times \sin \beta
\]
\[
1 + \sin \alpha = \frac{\cos \alpha}{\sin \beta} \times \cos \beta
\]
Multiplying the above two equations, we get
\[
1 - \sin^2 \alpha = \cos^2 \alpha
\]
or
\[
\cos^2 \alpha + \sin^2 \alpha = 1.
\]
As we expected.

\subsection*{Case 5:}

\begin{figure}[H]
    \centering
    \scalebox{1.2}{\input{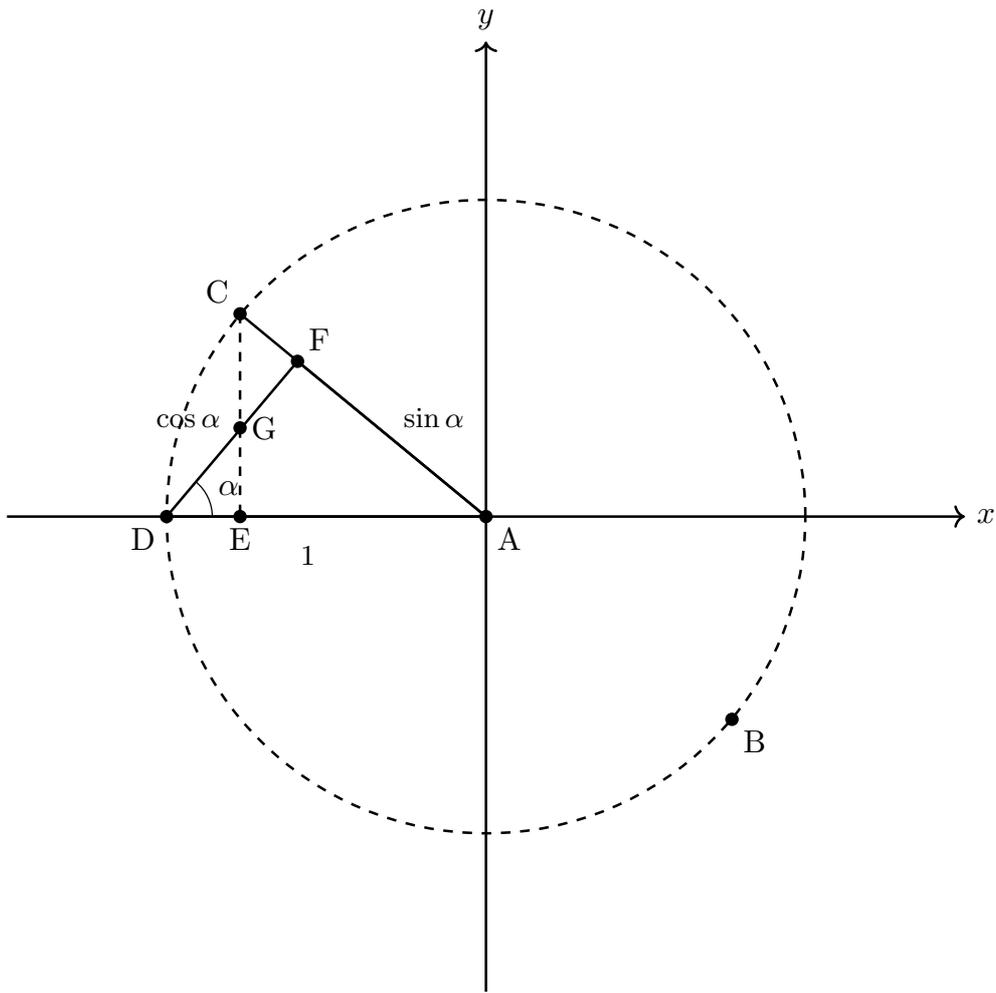}}
    \renewcommand{\thefigure}{A2.5a} % Temporarily redefine figure numbering
    \captionsetup{labelformat=empty} % Removes default label format
    \caption{Fig.~A2.5a: An alternative proof of the Pythagorean Theorem using the properties of a triangle inscribed in a circle}
    \label{Figure_A2.5a}
    \renewcommand{\thefigure}{\arabic{figure}} % Restore standard numbering
\end{figure}

In Fig.~\ref{Figure_A2.5a}, points \(C\) and \(D\) are on a circle with a radius of 1 unit. The Primary Gasing Triangle \(\triangle ADF\) has a hypotenuse \(AD = 1\) unit. 

Next, we draw line \(CE\) perpendicular to \(DA\) so that Primary Gasing Triangle \(\triangle ACE\) is similar to triangle \(\triangle ADF\). 

Since \(\angle ACE = \alpha\), then \(AE = \sin \alpha\). Thus, \(DE = 1 - \sin \alpha\). See Fig.~\ref{Figure_A2.5b}.

\begin{figure}[H]
    \centering
    \scalebox{1}{\input{figures/Figure_A2.5b}}
    \renewcommand{\thefigure}{A2.5b} % Temporarily redefine figure numbering
    \captionsetup{labelformat=empty} % Removes default label format
    \caption{Fig.~A2.5b: Side lengths of the triangles from Figure A2.5a}
    \label{Figure_A2.5b}
    \renewcommand{\thefigure}{\arabic{figure}} % Restore standard numbering
\end{figure}

Triangles \(\triangle ADF\) and \(\triangle DGE\) are similar. Using the scaling factor \(\frac{\sin \alpha}{\cos \alpha}\) on triangle \(\triangle ADF\), we obtain triangle \(\triangle DGE\). Since \(DE = 1 - \sin \alpha\), then
\[
EG = \frac{\sin \alpha}{\cos \alpha} \times (1 - \sin \alpha).
\]

Next, \(FC = 1 - \sin \alpha\), so
\[
GC = \frac{1}{\cos \alpha} \times (1 - \sin \alpha).
\]

Now, since \(CE = \cos \alpha\) and \(EC = EG + GC\), we have the equation
\[
\cos \alpha = \frac{\sin \alpha}{\cos \alpha} \times (1 - \sin \alpha) + \frac{1}{\cos \alpha} \times (1 - \sin \alpha)
\]
or
\[
\cos^2 \alpha = \sin \alpha\,(1 - \sin \alpha) + (1 - \sin \alpha).
\]
Simplifying, we get
\[
\cos^2 \alpha = 1 - \sin^2 \alpha.
\]
Thus, we have proven
\[
\cos^2 \alpha + \sin^2 \alpha = 1
\]
as we expected.

\subsection*{Case 6:}

\begin{figure}[H]
    \centering
    \scalebox{1.1}{\input{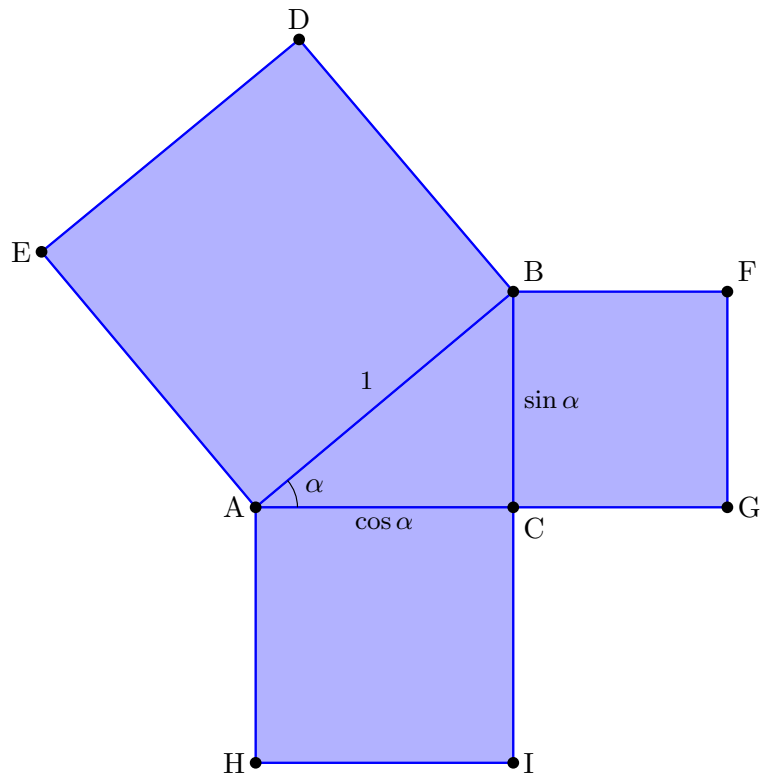}}
    \renewcommand{\thefigure}{A2.6a} % Temporarily redefine figure numbering
    \captionsetup{labelformat=empty} % Removes default label format
    \caption{Fig.~A2.6a: Three squares and one triangle for proving the Pythagorean Theorem}
    \label{Figure_A2.6a}
    \renewcommand{\thefigure}{\arabic{figure}} % Restore standard numbering
\end{figure}

In Fig.~\ref{Figure_A2.6a}, we have the Primary Gasing Triangle \(\triangle ABC\) with three squares \(\square HICA\), \(\square CGFB\), and \(\square ABDE\). 

Here, the hypotenuse of the Primary Gasing Triangle \(\triangle ABC\) is \(AB = 1\), the length of \(AC = \cos \alpha\), and the length of \(CB = \sin \alpha\). 

From the figure, we get the following:
\[
\text{The area of square } \square HICA = \cos^2 \alpha
\]
\[
\text{The area of square } \square CGFB = \sin^2 \alpha
\]
\[
\text{The area of square } \square ABDE = 1
\]
Now, we draw a line \(CK\) perpendicular to line \(AB\).

\begin{figure}[H]
    \centering
    \scalebox{1.1}{\input{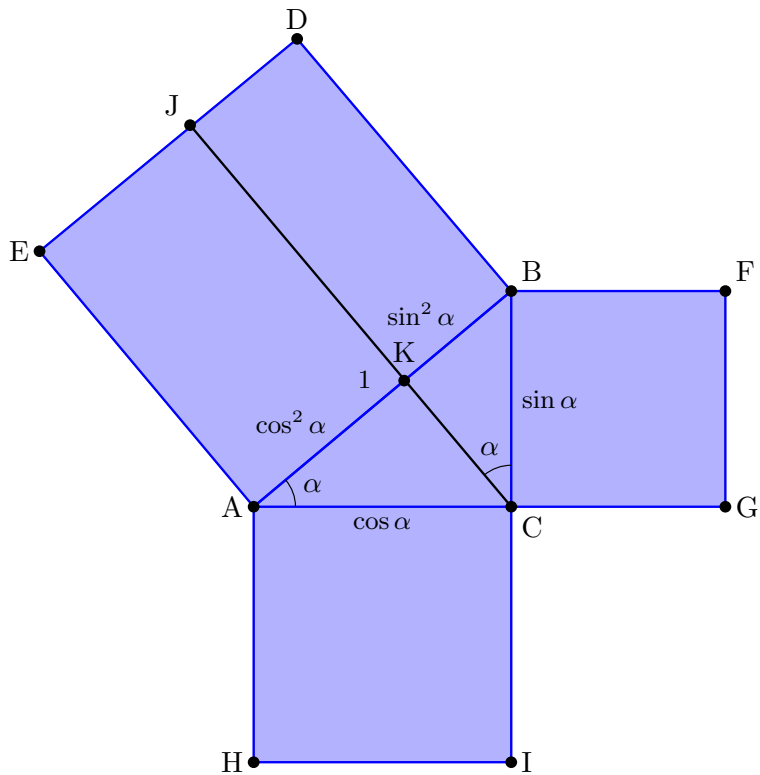}}
    \renewcommand{\thefigure}{A2.6b} % Temporarily redefine figure numbering
    \captionsetup{labelformat=empty} % Removes default label format
    \caption{Fig.~A2.6b: Side lengths from Figure A2.6a for proving the Pythagorean Theorem}
    \label{Figure_A2.6b}
    \renewcommand{\thefigure}{\arabic{figure}} % Restore standard numbering
\end{figure}

Triangle \(\triangle CBK\) is similar to the Primary Gasing Triangle \(\triangle ABC\). 

Since side \(BC = \sin \alpha\), the length of side \(BK = \sin^2 \alpha\) (using the scaling factor \(\sin \alpha\)). 

Similarly, since triangles \(\triangle ABC\) and \(\triangle ACK\) are similar, the length of side \(AK = \cos^2 \alpha\) (using the scaling factor \(\cos \alpha\)).

Since the length of side \(AE = BD = 1\), then:
\[
\text{The area of } AKJE = AK \times AE = \cos^2 \alpha \times 1 = \cos^2 \alpha
\]
\[
\text{The area of } KBDJ = KB \times BD = \sin^2 \alpha \times 1 = \sin^2 \alpha
\]
Since the area of \(\square ABDE = \text{the area of } AKJE + \text{the area of } KBDJ\), thus:
\[
1 = \cos^2 \alpha + \sin^2 \alpha
\]

This equation simultaneously proves the Pythagorean theorem and shows that
\[
\text{The area of } \square ABDE = \text{The area of } \square HICA + \text{The area of } \square CGFB.
\]
The latter is often used for a geometric proof of the Pythagorean theorem.

\subsection*{Case 7:}

\begin{figure}[H]
    \centering
    \scalebox{1.2}{\input{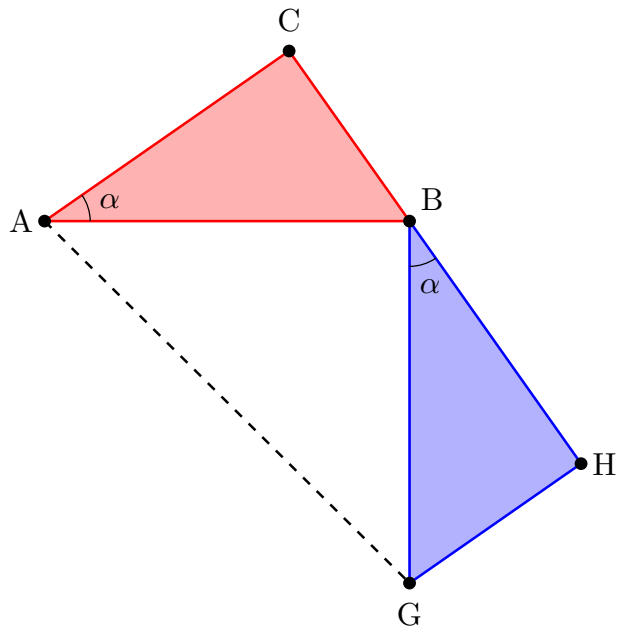}}
    \renewcommand{\thefigure}{A2.7a} % Temporarily redefine figure numbering
    \captionsetup{labelformat=empty} % Removes default label format
    \caption{Fig.~A2.7a: Three triangles for proving the Pythagorean Theorem}
    \label{Figure_A2.7a}
    \renewcommand{\thefigure}{\arabic{figure}} % Restore standard numbering
\end{figure}

In Fig.~\ref{Figure_A2.7a}, we have a trapezium \(ACHG\) divided into three right triangles: Primary Gasing Triangle \(\triangle ABC\), Primary Gasing Triangle \(\triangle BGH\), and triangle \(\triangle AGB\). The Primary Gasing Triangles \(\triangle ABC\) and \(\triangle BGH\) are congruent, with side \(AB = 1\) and angle \(\angle CAB = \alpha\).

\begin{figure}[H]
    \centering
    \scalebox{1.2}{\input{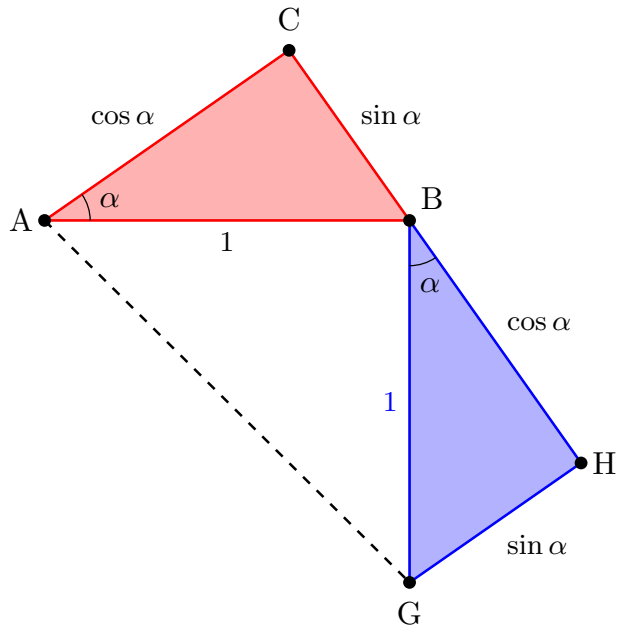}}
    \renewcommand{\thefigure}{A2.7b} % Temporarily redefine figure numbering
    \captionsetup{labelformat=empty} % Removes default label format
    \caption{Fig.~A2.7b: Side lengths of the triangles from Figure A2.7a for proving the Pythagorean Theorem.}
    \label{Figure_A2.7b}
    \renewcommand{\thefigure}{\arabic{figure}} % Restore standard numbering
\end{figure}

Since \(\triangle ABC\) is a Primary Gasing Triangle, the lengths of the sides are \(BC = \sin \alpha\) and \(AC = \cos \alpha\). Because triangles \(\triangle ABC\) and \(\triangle BGH\) are congruent, \(\angle GBH = \alpha\) (since line \(AB\) is perpendicular to \(BG\)). Therefore, \(BH = \cos \alpha\) and \(GH = \sin \alpha\), as shown in Fig.~\ref{Figure_A2.7b}.

We see that the area of trapezoid \(ACHG = 2 \times \text{The area of } \triangle ABC + \text{The area of } \triangle AGB\), or
\[
\frac{(\sin \alpha + \cos \alpha)(\sin \alpha + \cos \alpha)}{2} = 2 \times \frac{(\sin \alpha \cos \alpha)}{2} + \frac{1}{2}.
\]
Expanding this, we get
\[
\sin^2 \alpha + 2 \sin \alpha \cos \alpha + \cos^2 \alpha = 2 \sin \alpha \cos \alpha + 1
\]
or
\[
\sin^2 \alpha + \cos^2 \alpha = 1
\]
as we expected.

\subsection*{Case 8:}

Here, we will work on the proof of the Pythagorean Theorem as performed by Calcea Johnson and Ne'Kiya Jackson \autocite{sloman_2023} using Primary Gasing Triangles.

\begin{figure}[H]
    \centering
    \scalebox{0.7}{\input{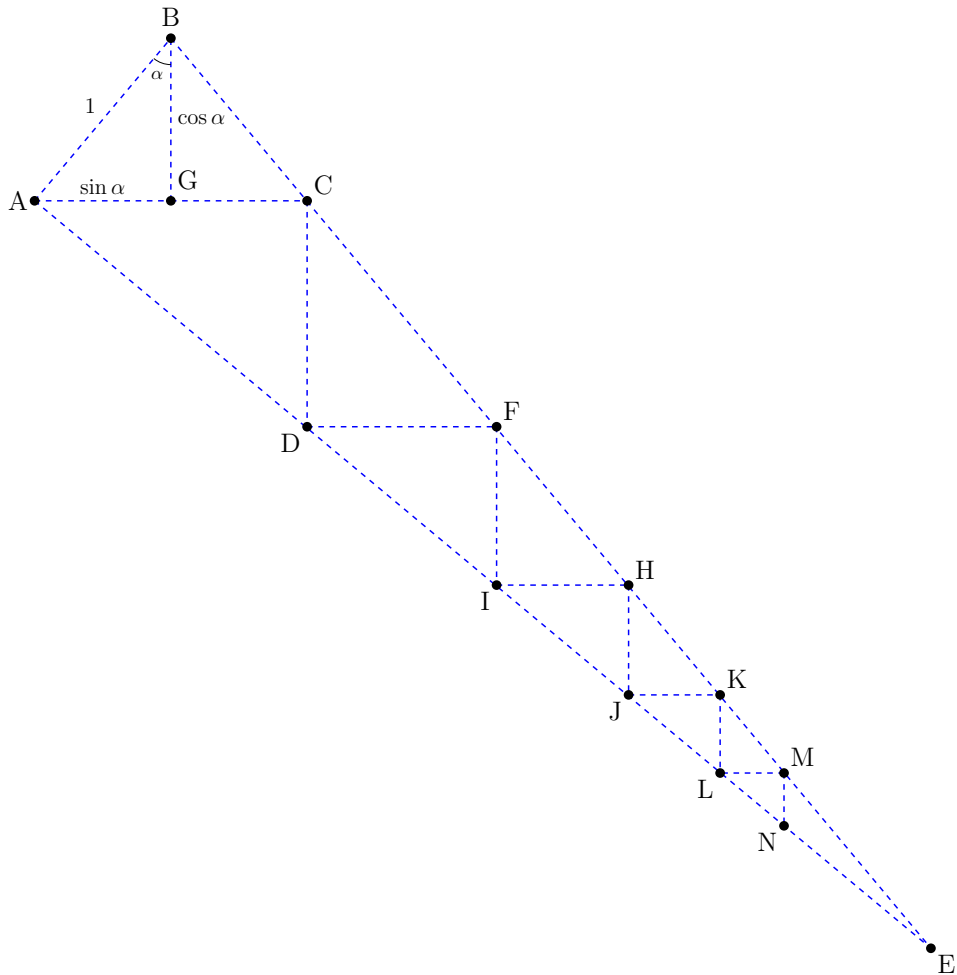}}
    \renewcommand{\thefigure}{A2.8a} % Temporarily redefine figure numbering
    \captionsetup{labelformat=empty} % Removes default label format
    \caption{Fig.~A2.8a: Calcea Johnson and Ne'Kiya Jackson's method for proving the Pythagorean Theorem, adapted using the Gasing Triangle}
    \label{Figure_A2.8a}
    \renewcommand{\thefigure}{\arabic{figure}} % Restore standard numbering
\end{figure}

Here, we take the Gasing Triangle \(\triangle BAG\) with \(BA = 1\) unit, \(GA = \sin \alpha\), and \(BG = \cos \alpha\). 

Triangles \(\triangle BCG\) and \(\triangle BAG\) are congruent, so \(AC = 2 \sin \alpha\).

Triangles \(\triangle BAG\) and \(\triangle DAC\) are similar. We use the scaling factor \(\frac{\sin \alpha}{\cos \alpha}\) on triangle \(\triangle BAG\) to obtain triangle \(\triangle DAC\). With this scale factor, we get the lengths \(AD = \frac{2 \sin \alpha}{\cos \alpha}\) and \(CD = \frac{2 \sin \alpha}{\cos \alpha} \times \sin \alpha\).

Next, since \(\angle DCF = \alpha\), and we already know the length of \(CD\), by using the similarity between triangles \(\triangle DAC\) and \(\triangle FCD\), we get the lengths
\[
FC = \frac{2 \sin^2 \alpha}{\cos^2 \alpha} \times \sin \alpha
\]
and
\[
DF = \frac{2 \sin^2 \alpha}{\cos^2 \alpha} \times \sin \alpha.
\]

We can continue this process until point \(E\). The result is as shown in the following figure.

\begin{figure}[H]
    \centering
    \scalebox{0.7}{\input{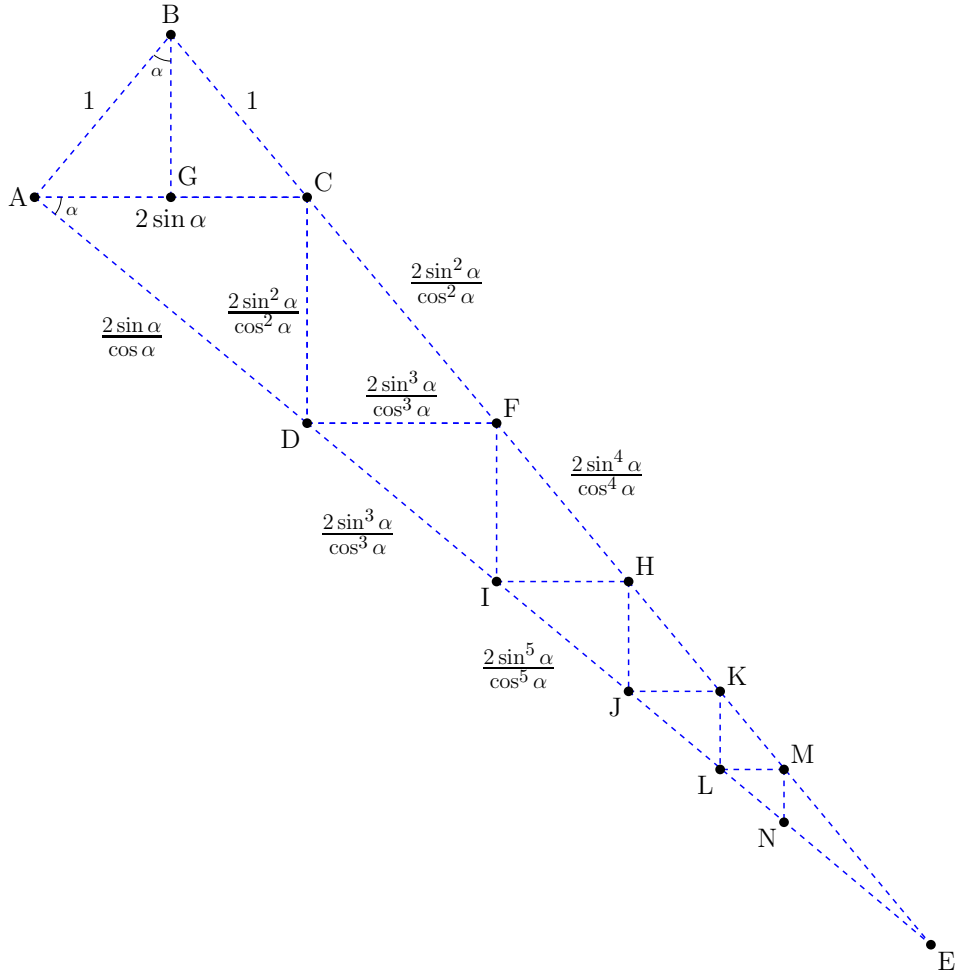}}
    \renewcommand{\thefigure}{A2.8b} % Temporarily redefine figure numbering
    \captionsetup{labelformat=empty} % Removes default label format
    \caption{Fig.~A2.8b: Side lengths of the triangles from Figure A2.8a}
    \label{Figure_A2.8b}
    \renewcommand{\thefigure}{\arabic{figure}} % Restore standard numbering
\end{figure}

The length of \(EA\) is calculated as:
\begin{align}
    EA &= \frac{2 \sin \alpha}{\cos \alpha} 
    \left(1 + \frac{\sin^2 \alpha}{\cos^2 \alpha} 
    + \frac{\sin^4 \alpha}{\cos^4 \alpha} + \cdots \right) \notag \\
    &= \frac{2 \sin \alpha}{\cos \alpha} 
    \cdot \frac{1}{1 - \frac{\sin^2 \alpha}{\cos^2 \alpha}} \notag \\
    &= \frac{2 \sin \alpha \cos \alpha}{\cos^2 \alpha - \sin^2 \alpha} \notag
\end{align}

Similarly, the length of \(EB\) is given by:
\begin{align}
    EB &= 1 + \frac{2 \sin^2 \alpha}{\cos^2 \alpha} 
    \left(1 + \frac{\sin^2 \alpha}{\cos^2 \alpha} 
    + \frac{\sin^4 \alpha}{\cos^4 \alpha} + \cdots \right) \notag \\
    &= 1 + \frac{2 \sin^2 \alpha}{\cos^2 \alpha} 
    \cdot \frac{1}{1 - \frac{\sin^2 \alpha}{\cos^2 \alpha}} \notag \\
    &= 1 + \frac{2 \sin^2 \alpha}{\cos^2 \alpha - \sin^2 \alpha} \notag \\
    &= \frac{\cos^2 \alpha + \sin^2 \alpha}{\cos^2 \alpha - \sin^2 \alpha} \notag
\end{align}

\begin{figure}[H]
    \centering
    \scalebox{0.7}{\input{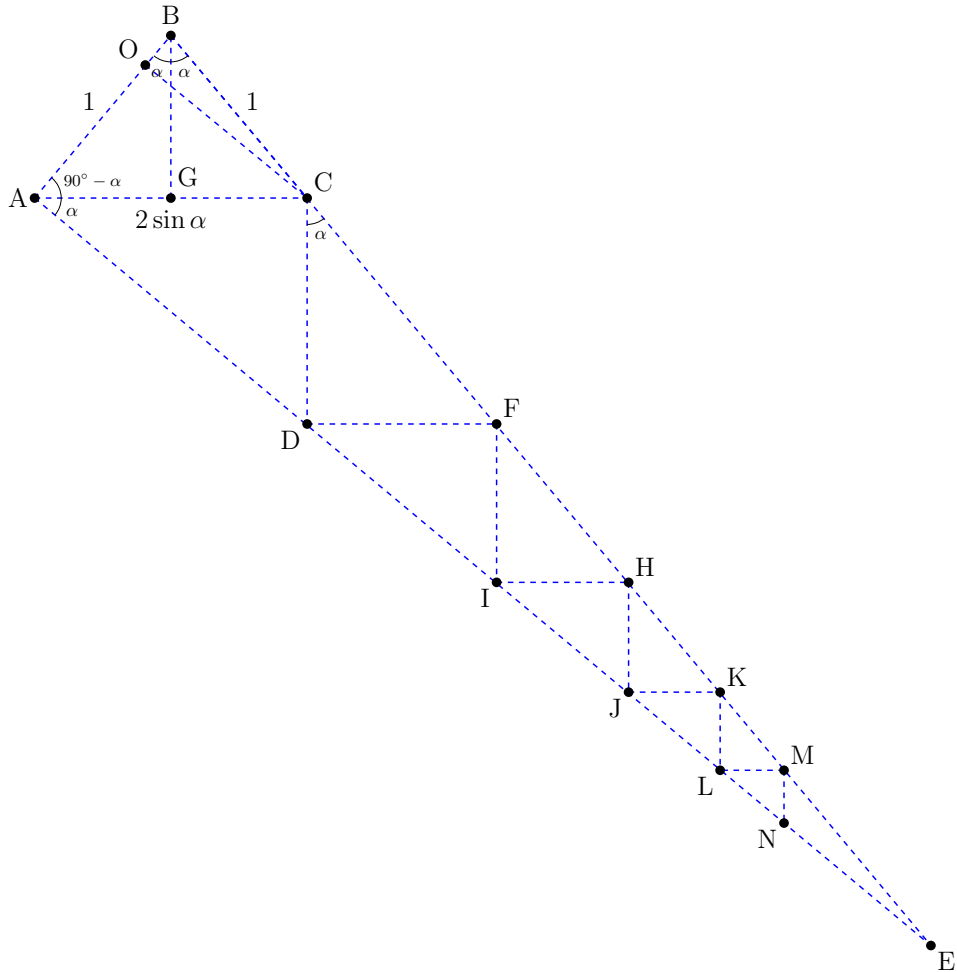}}
    \renewcommand{\thefigure}{A2.8c} % Temporarily redefine figure numbering
    \captionsetup{labelformat=empty} % Removes default label format
    \caption{Fig.~A2.8c: Further steps in the proof of the Pythagorean Theorem}
    \label{Figure_A2.8c}
    \renewcommand{\thefigure}{\arabic{figure}} % Restore standard numbering
\end{figure}

From the Gasing Triangle \(\triangle CBO\) (right-angled at \(O\)) with \(\angle CBO = 2\alpha\), we get that the length of \(OC = \sin 2\alpha\). 

From triangle \(\triangle CAO\) (right-angled at \(O\)) with length \(CA = 2 \sin \alpha\), we get the length of $OC = 2 \sin \alpha \sin\,(90^\circ - \alpha) = 2 \sin \alpha \cos \alpha$.

So here we get \(\sin 2\alpha = 2 \sin \alpha \cos \alpha\) (without using the sum identity of sine).

From triangle \( \triangle EBA \) with angle \( \angle EBA = \alpha \), we get
\[
EA = EB \sin 2\alpha = EB \times 2 \sin \alpha \cos \alpha.
\]

Substituting the values of \(EA\) and \(EB\), we have
\[
\frac{2 \sin \alpha \cos \alpha}{\cos^2 \alpha - \sin^2 \alpha}
= \frac{\cos^2 \alpha + \sin^2 \alpha}{\cos^2 \alpha - \sin^2 \alpha} \times 2 \sin \alpha \cos \alpha.
\]

Simplifying, we get
\[
1 = \cos^2 \alpha + \sin^2 \alpha
\]

as we expected.

\break
\printbibliography

@Article{de_villiers_2023,
  author    = {Michael de Villiers},
  title     = {Is a Trigonometric Proof Possible for the Theorem of Pythagoras?},
  journal   = {Learning and Teaching Mathematics},
  number    = {34},
  year      = {2023},
  pages     = {22-27},
}

@Article{shirali_2023,
  author    = {Shailesh Shirali},
  title     = {Two New Proofs of the Pythagorean Theorem},
  journal   = {At Right Angles},
  publisher = {Azim Premji University},
  year      = {2023},
}

@Article{sloman_2023,
  author    = {Sloman, Leila},
  title     = {2 High School Students Prove Pythagorean Theorem. Here’s What That Means},
  journal   = {Scientific American},
  year      = {2023},
  url       = {https://www.scientificamerican.com/article/2-high-school-students-prove-pythagorean-theorem-heres-what-that-means/},
}

@Book{bylen_ziegler_barnett_2011,
  author    = {K. E. Bylen and M.R. Ziegler and R.A. Barnett},
  title     = {Analytic Trigonometry with Applications, 11th Edition},
  publisher = {Wiley},
  year      = {2011},
}

@Book{gelfand_saul_2001,
  author    = {I.M. Gelfand and Mark Saul},
  title     = {Trigonometry},
  publisher = {Birkhäuser},
  address   = {Boston},
  year      = {2001},
}

@Book{gullberg_mathematics_1996,
  author    = {Gullberg, Jan},
  title     = {Mathematics from the Birth of Numbers},
  publisher = {W. W. Norton and Company, Inc.},
  address   = {New York},
  year      = {1996},
}

@Book{hall_knight_1928,
  author    = {H.S. Hall and S.R. Knight},
  title     = {Elementary Trigonometry},
  publisher = {The Macmillan Company of Canada Ltd.},
  address   = {Toronto},
  year      = {1928},
}

@Book{loomis_1968,
  author    = {Elisha Scott Loomis},
  title     = {The Pythagorean Proposition},
  publisher = {NCTM},
  year      = {1968},
}

@Book{simmons_1987,
  author    = {Simmons, George F.},
  title     = {Precalculus Mathematics in a Nutshell},
  publisher = {Barnes and Noble Books},
  address   = {New York},
  year      = {1987},
}

@Book{stewart_redlin_watson_2016,
  author    = {J. Stewart and L. Redlin and S. Watson},
  title     = {Precalculus: Mathematics for Calculus, 7th Ed.},
  publisher = {Cengage Learning},
  address   = {Australia},
  year      = {2016},
}

@Book{stitz_zeager_2012,
  author    = {Carl Stitz and Jeff Zeager},
  title     = {College Trigonometry},
  year      = {2012},
  note      = {Lakeland Community College and Lorain County Community College, January 8, 2012},
}

@Book{sullivan_2016,
  author    = {M. Sullivan},
  title     = {Trigonometry: A Unit Circle Approach, 10th Ed.},
  publisher = {Pearson},
  address   = {Boston},
  year      = {2016},
}

@Book{sundstrom_schlicker_2019,
  author    = {Ted Sundstrom and Steven Schlicker},
  title     = {Trigonometry},
  publisher = {Creative Commons},
  address   = {California},
  year      = {2019},
}

@Book{van_brummelen_2009,
  author    = {Van Brummelen, Glen},
  title     = {The Mathematics of The Heavens and The Earth: The Early History of Trigonometry},
  publisher = {Princeton University Press},
  address   = {New Jersey},
  year      = {2009},
}

@Online{gasing_academy,
  author    = {{Gasing Academy}},
  title     = {Gasing Academy},
  year      = {n.d.},
  url       = {https://www.gasingacademy.org},
}

@PhdThesis{van_sickle_2011,
  author       = {Van Sickle, Jenna},
  title        = {A History of Trigonometry Education in the United States: 1776-1900},
  year         = {2011},
  school       = {Columbia University},
  type         = {Ph.D. dissertation},
  publisher    = {ProQuest Dissertations \& Theses Global},
  note         = {Publication No. 3453860}
}

\end{document}